\documentclass[12pt]{amsart}

\usepackage{amssymb,amsfonts,amsmath,amsopn,amstext,amscd,color,latexsym,
mathrsfs,mathtools,stackrel,url,verbatim}
\usepackage{mathabx}
\usepackage{tikz-cd}
\usepackage[utf8]{inputenc}
\usepackage[all, cmtip]{xy}
\usepackage[pagewise,displaymath,mathlines]{lineno}
\usepackage[active]{srcltx}
\usepackage{xcolor}
\usepackage[normalem]{ulem}  

\usepackage[pdftex, colorlinks=true, linktoc=all]{hyperref}

\usepackage{cleveref}

\newcommand{\nc}{\newcommand}
\nc{\rc}{\renewcommand}

\nc{\N}{\mathbb{N}}
\nc{\C}{\mathbb{C}}
\nc{\Z}{\mathbb{Z}}
\nc{\wt}{\widetilde}
\nc{\wh}{\widehat}
\nc{\Q}{\mathbb{Q}}
\nc{\Qbar}{\overline{\mathbb{Q}}}
\nc{\R}{\mathbb{R}}
\nc{\Fp}{\mathbb{F}_p}
\nc{\Fpbar}{\overline{\mathbb{F}_p}}
\nc{\units}{^{\times}}
\nc{\inverse}{^{-1}}
\nc{\m}{\mathbf{m}}
\rc{\t}{\mathrm}
\nc{\s}{\mathcal}
\nc{\Der}{\t{Der}}
\nc{\Hom}{\t{Hom}}
\nc{\End}{\t{End}}
\nc{\Aut}{\t{Aut}}
\nc{\GL}{\t{GL}}
\nc{\SL}{\t{SL}}
\nc{\Mn}{\t{M}_n}
\nc{\Mg}{\t{M}_g}
\nc{\Gal}{\t{Gal}}
\nc{\absgal}{\Gal(\bar{\Q}/\Q)}
\nc{\gr}{\t{gr}}
\rc{\P}{\mathbb{P}}
\rc{\div}{\t{div}} 
\nc{\ord}{\t{ord}}
\nc{\sm}{\t{C}^{\infty}}
\nc{\pd}{\partial}
\nc{\df}{C^{\infty}}
\nc{\vt}{\vartheta}
\rc{\O}{\s{O}}
\nc{\HH}{\mathbb{H}}
\nc{\V}{\t{V}}
\nc{\D}{\t{D}}
\nc{\Kn}{K[T_1, T_2,...,T_n]}
\nc{\Cn}{\C[T_1,T_2,...,T_n]}
\nc{\rad}{\t{rad}}
\nc{\An}{A[T_1,T_2,...,T_n]}
\nc{\Mor}{\t{Mor}}
\nc{\alg}{\t{alg}}
\nc{\Sp}{\t{Spec}}
\nc{\spr}{\t{Sp}_{red}}
\nc{\AL}{\t{AL}}
\nc{\rng}{\t{Rings}}
\nc{\MA}{\Mor_{K-\alg}}
\rc{\MR}{\Mor_{\t{Rings}}}
\nc{\A}{\mathbb{A}}
\nc{\AXY}{A[X_1,...,X_n,Y_1,...,Y_m]}
\nc{\omba}{\Omega_{B/A}}
\nc{\Fq}{\mathbb{F}_q}
\nc{\Ga}{\mathbb{G}_a}
\nc{\Di}{Dieudonne}
\nc{\h}{\mathfrak{h}}
\nc{\Res}{\t{Res}}
\nc{\Rd}{\t{R}}
\nc{\Div}{\t{Div}}
\nc{\Cl}{\t{Cl}}
\nc{\g}{\gamma}
\rc{\dot}[2]{\langle#1,#2\rangle}
\nc{\ufr}[1]{\lceil#1\rceil}
\rc{\l}{\mathit{l}}
\nc{\tr}{\mathit{t}}
\nc{\G}{\Gamma}
\rc{\Re}{\t{Re}}
\rc{\Im}{\t{Im}}
\nc{\Ell}{\t{Ell}}
\nc{\Pic}{\t{Pic}}
\nc{\p}[1]{\textbf{Problem}\,#1:}
\nc{\sq}[1]{#1\units/(#1\units)^2 }
\nc{\lr}{\leftarrow}

\nc{\Gm}{\mathbb{G}_m}
\nc{\bsl}{\backslash}
\rc{\l}{\mathit{l}}
\rc{\ll}{L_{\lambda}}
\nc{\la}{\lambda}
\nc{\Heis}{\t{Heis}}

\rc{\H}{\t{H}}
\rc{\bar}{\overline}
\nc{\xr}{\xrightarrow}
\nc{\xl}{\xleftarrow}

\nc{\hkr}{\hookrightarrow}

\newtheorem{thm}{Theorem}[subsection]  
\newtheorem{lem}[thm]{Lemma}
\newtheorem{prop}[thm]{Proposition}

\newtheorem{cor}[thm]{Corollary}

\theoremstyle{definition}
\newtheorem{ex}[thm]{Example}

\newtheorem{rmk}[thm]{Remark}

\newtheorem{dfn}[thm]{Definition}
\newtheorem{pp}[thm]{Problem}
\newtheorem{question}{Question}[section]

\nc{\bp}{\begin{proof}}
\nc{\ep}{\end{proof}}
\nc{\bprop}{\begin{prop}}
\nc{\eprop}{\end{prop}}
\nc{\br}{\begin{rmk}}
\nc{\er}{\end{rmk}}
\nc{\bl}{\begin{lem}}
\nc{\el}{\end{lem}}
\nc{\bd}{\begin{dfn}}
\nc{\ed}{\end{dfn}}
\nc{\bcor}{\begin{cor}}
\nc{\ecor}{\end{cor}}
\nc{\be}{\begin{enumerate}}
\nc{\ee}{\end{enumerate}}
\nc{\bt}{\begin{thm}}
\nc{\et}{\end{thm}}
\nc{\bpp}{\begin{pp}}
\nc{\epp}{\end{pp}}
\nc{\bex}{\begin{ex}}
\nc{\eex}{\end{ex}}
\nc{\beq}{\begin{equation}}
\nc{\eeq}{\end{equation}}
\nc{\disk}{\Delta}
\nc{\pdisk}{\Delta^*}

\begin{document}
\title{Local Monodromy of Constructible Sheaves.}
\author{Madhav V. Nori}
\address{Department of Mathematics, University of Chicago, 5734 University Ave., 
Chicago, IL-60637, U.S.A.}
\email{nori@math.uchicago.edu}
\author{Deepam Patel}
\address{Department of Mathematics, Purdue University,
150 N. University Street, West Lafayette, IN 47907, U.S.A.}
\email{patel471@purdue.edu}

\thanks{D. Patel was partially supported by a Simons Collaboration Grant.}

\begin{abstract} 
Given a morphism $f: X \rightarrow S$ of complex algebraic varieties and a constructible sheaf $\s{G}$ on $X$, we compute the local monodromy of $Rf_*(\s{G})$ and $Rf_!(\s{G})$ in terms of the local monodromy of $\s{G}$. Our results generalize previous results by Brieskorn, Borel, Clemens, Deligne, Landman, Griffiths, Grothendieck, and Kashiwara in the setting of quasi-unipotent sheaves. In the following, we consider the general setting of sheaves of $R$-modules for a commutative noetherian ring $R$, and give applications to computing local monodromy of abelian covers in a {\it uniform manner}. We also obtain applications in the context of `generalized Alexander modules' and intersection cohomology with torsion coefficients.

\end{abstract}
\maketitle
\tableofcontents

\section{Introduction}

In the following article, we investigate the behavior of eigenvalues arising from the action of {\it analytic loops} on constructible sheaves under Grothendieck's six-functor formalism. We work in the setting of constructible sheaves of $R$-modules for a commutative noetherian ring $R$ of finite homological dimension. In particular, we study situations arising from {\it non-geometric local systems}. An interesting feature of our applications is that the consideration of local monodromy for possibly {\it non-geometric} local systems and with coefficients in $R$-modules has implications for certain geometric questions. More precisely, consider a commutative diagram of complex algebraic varieties:
$$
\xymatrix{
X \ar[r]^{F} \ar[dr]^{f} & G \ar[d]^{\pi} \\
      &   S}
$$
with $S$ a smooth curve, and $G$ a semi-abelian scheme (over $S$). Let $n_G: G \rightarrow G$ be the multiplication by $n$ map, and consider the resulting etale covers $n_X: X_n \rightarrow X$ obtained via base change along $F$. Let $f_n := f \circ n_X$, $s_0 \in S$, and let $\Delta \rightarrow S$ be a small disk centered at $s_0$. Up to further shrinking of the disk, we may assume that $R^if_{n,*}\Z$ is a local system when restricted to the punctured disk. This data gives a sequence of local monodromy representations:
$$
\rho_{n}: \pi_1(\Delta^*) \rightarrow {\rm GL}((R^if_{n,*}\Z)_t)
$$
for a general point $t \in \Delta^*$. Then Grothendieck's {\it Local Monodromy Theorem} (\cite{SGA7-V1}) states that the eigenvalues of the local monodromy action on $R^if_{n,*}\Z$ are roots of unity. This leads to the following natural question:

\begin{question}\label{question:questionintro}
With notation as above, which roots of unity appear? Can these be obtained in a uniform manner (as $n$-varies)?
\end{question}

We consider such examples in Section \ref{sec:alexmodules} and show that these roots can be computed explicitly and uniformly in $n$ as an application of our general results. We also obtain applications to the local monodromy of Alexander modules and to the local monodromy action on intersection cohomology (with coefficients in arbitrary fields). Our results generalize various results in the literature on quasi-unipotence of local monodromy. We discuss these in detail below.

\subsection{Main Results}
Let $X$ be a variety over $\C$, and $\s{G}$ be a constructible sheaf of $K$-vector spaces where $K$ is an algebraically closed field. In the following, we let $\s{O}^{an} := \C\{t\}$ denote the ring of convergent power series in the variable $t$, and $F^{an}$ denote its fraction field. 

\bd\label{defn:analyticloop}
\begin{enumerate}
\item An {\it analytic loop} of $X$ is a morphism of $\C$-schemes 
$$\gamma: \Sp(F^{an}) \rightarrow X.$$ We denote by
${\rm AL}(X)$ the set of analytic loops of $X$. In the following, we shall sometimes refer to analytic loops as loops. 
\item Let $[n]: \Sp(F^{an}) \rightarrow \Sp(F^{an})$ induced by $t \to t^n$. For $\g \in {\rm AL}(X)$ we denote by $\g^n := \g \circ [n]$.
\end{enumerate}
\ed

The loop $\gamma \in \AL(X)$ gives an analytic map $h: \pdisk \rightarrow X$ from a small punctured disk, and (up to shrinking the disk) we may assume that $h^*\s{G}$ is locally constant.\footnote{We refer the reader to section \ref{sec:bmon}, before definition \ref{defn:curveloops} for an explanation of these facts.} After fixing base points, and considering the canonical generator (i.e. counterclockwise loop around the origin) $T \in \pi_1(\pdisk)$, we may consider the set of eigenvalues of the resulting monodromy action on the (stalk of the) given local system. The set of eigenvalues thus obtained is independent of the choice of disk or base point, and we let $\spr(\gamma, \s{G}) \subset K^{\times}$ (the `{\it reduced spectrum}' of the given loop) denote the resulting set of eigenvalues. We set ${\rm BSp}(\s{G}) : = \bigcup_{\gamma \in \AL(X)} \spr(\g,\s{G})$ (the `{\it boundary spectrum}' of $\s{G})$.

\br
\begin{enumerate}\label{rem:bmonprop}
   \item \label{rem:bmonprop1} If $\s{G}$ is non-zero, then $1 \in {\rm BSp}(\s{G})$. If $\s{G}$ is a non-zero constant local system, then ${\rm BSp}(\s{G}) = \{1\}.$ 
   \item \label{rem:bmonprop2 } Moreover, if $c \in {\rm BSp}(\s{G})$, then $c^k \in {\rm BSp}(\s{G}) $ for all $k >0$.
\end{enumerate}
\er

The main result of this article is the following theorem on the behavior of reduced spectra (and the boundary spectrum) under push-forwards. In the following, given a subset $M \subset K^{\times}$ and an integer $r > 0$, we set $M^{\frac{1}{r}} := \{ c \in K^{\times} | c^r \in M\}$ and let $M^{+}$ denote the monoid generated by the elements of $M \subset K^{\times}$.

\bt \label{thm:maintthm}
Let $f: X \rightarrow S$ be a morphism of a complex algebraic varieties, and $\s{G}$ a constructible sheaf of $K$-vector spaces on $X$.

\begin{enumerate}

\item\label{thm:maintthm1} Let $\gamma \in {\rm AL}(S)$. Then there is an integer $r >0 $ and a finite subset $M \subset {\rm AL}(X)$ (both depending only on $\s{G}$, $f$ and $\gamma$), such that 

$$\spr(\gamma, R^qf_?(\s{G})) \subset \{\lambda | \lambda^r \in \bigcup_{\gamma' \in M} \spr(\gamma',\s{G})\},$$

where $? \in \{*,!\}$.

\item\label{thm:maintthm2} Suppose $\dim(S) =1$. There is an integer $r > 0$ (depending on $\s{G}$ and $f$), such that for all $q$, ${\rm BSp}(R^qf_?\s{G}) \subset {\rm BSp}(\s{G})^{\frac{1}{r}}$, where $? \in \{*,!\}$.
\item\label{thm:maintthm3} In general (for $\dim(S) \geq 1$), there is an integer $r' >0 $ (depending on $\s{G}$ and $f$), such that for all $q$,  ${\rm BSp}(R^qf_?\s{G})^{+} \subset ({\rm BSp}(\s{G})^{+})^{\frac{1}{r'}}$, where $? \in \{*,!\}$.
\end{enumerate}
\et

Note that, as a consequence of the fact that $R^qf_!(\s{G})$ and $R^qf_*(\s{G})$ vanish for large $q$, it is enough to show that such an integer exists for a fixed $q$. Similarly, it follows that there exists an integer $r$ that will work for both $f_!$ and $f_*$.

\br
In the text, we shall work more generally with (bounded) constructible complexes, i.e. objects of the derived category. The definitions and results above generalize immediately to this setting by passing to the corresponding cohomology sheaves. 
\er

The aforementioned theorem generalizes the classical theorem on quasi-unipotence of local monodromy due to Grothendieck (\cite{SGA7-V1}), and related generalizations. We shall briefly recall the history and related results in section \ref{sec:history/related} below.\\

Our second result is an analogous assertion in the setting of constructible sheaves of $R$-modules where $R$ is a commutative noetherian ring of finite homological dimension. In order to state the result, we first introduce some terminology. Let $X$ be as before and $\s{G}$ be a constructible sheaf of $R$-modules (on $X$). Given an analytic loop $\gamma \in {\rm AL}(X)$, one has an induced $R$-linear map 
$T: h^*\s{G}_t \rightarrow h^*\s{G}_t$ (with $h$
 as in the line following Definition \ref{defn:analyticloop}, and $t \in \pdisk$). In particular, we may view $h^*\s{G}_t$ as an $R[x]$-module. We note that (up to isomorphism) this module is independent of the choice of $h$ or base points (and $t$).

 \bd\label{defn:spRmod}
 Let $X$ be a scheme of finite type over $\C$, $\s{G}$ a constructible sheaf of $R$-modules, $\g \in {\rm AL}(X)$, and $h: \pdisk \rightarrow X$ the map associated to $\g$ so that $h^*(\s{G})$ is a local system. 
 \begin{enumerate}
 \item With notation as above, we set ${\rm Sp}(\gamma,\s{G}) \subset \Sp{(R[x])}$ to be the scheme theoretic support of $h^*\s{G}_t$ ($t \in \pdisk$) and by $I(\gamma,\s{G}) \subset R[x]$ the corresponding ideal.\footnote{This notation is consistent with the previous notation for the reduced spectrum. Specifically, if $R =K$ is an algebraically closed field, then the set of closed points of the underlying reduced scheme of the spectrum defined here is the reduced spectrum defined previously.}
 \item Given a positive integer $r$, and a closed subscheme $W \subset \A^1_R$, let $W^{[1/r]}$ denote the scheme theoretic inverse image of $W$ under the morphism $ \A^1_R \rightarrow \A^1_R$ induced by $x \to x^r$.
 \end{enumerate}
\ed

 \br\label{rem:sumofsubschemes}
 In the following, we consider {\it sums} of subschemes. Given a scheme $X$ and a (finite) collection of closed subschemes $Z_{\alpha} \subset X$, we {\bf define} the {\it sum} $\sum_{\alpha} Z_{\alpha}$ to be the closed subscheme defined by the product of ideals $I_{\alpha}$ defining each $Z_{\alpha}$. 
\er

\bt\label{thm:mainthm2}
Let $f: X \rightarrow S$ be a morphism of complex algebraic varieties, and $\s{G}$ a constructible sheaf of $R$-modules.
Let $\gamma \in {\rm AL}(S)$. Then there is a finite set (denoted by $M$) of pairs $(\g',n_{\g'})$ where $\g' \in {\rm AL}(X)$ and $n_{\g'}$ is a positive integer (depending only on $\s{G}$, $f$ and $\gamma$), such that 

$${\rm Sp}(\gamma, R^qf_?(\s{G})) \subset  \sum_{(\gamma',n_{\g'}) \in M} {\rm Sp}(\gamma',\s{G})^{[1/n_{\gamma'}]},$$

where $? \in \{*,!\}$.

\et

If $R =K$ is an algebraically closed field, then Theorem \ref{thm:mainthm2} implies Theorem \ref{thm:maintthm} (1) by taking the underlying reduced subscheme. On the other hand, the statement at the level of subschemes is {\it stronger} even in this case, since it controls the level of quasi-unipotency.

\br
For a natural number $r$ and a constructible sheaf $\s{G}$ of $R$-modules on $X$, one may also define a version of ${\rm BSp}(\s{G})^{\frac{1}{r}}$ as follows. We define ${\rm BSp}_{R}(\s{G})^{\frac{1}{r}}$ to be the smallest collection $M$ of closed subschemes of $\A^1_R$ with the following properties:
\\(I) If  $B\in M$ and $A$ is a closed subscheme of $B$, then $A\in M$
\\(II) If $A_1,...,A_m\in M$, then $\Sigma_{i=1}^{i=m}A_i\in M$ 
\\(III) ${\rm Sp}(\g,\s{G})^{[1/r]}\in M$   for every $\g\in {\rm AL}(X)$
\\We define ${\rm BSp}_R(\s{G})$ to be ${\rm BSp}_R(\s{G})^{\frac{1}{r}}$ where $r=1$. With this definition, the analog of Theorem \ref{thm:maintthm} (2) holds for sheaves of $R$-modules. 
\er

\subsection{Applications}
We survey some immediate applications of our main results. 

\subsubsection{Application to Integral transforms and Intersection cohomology}
As an immediate application of our main results, we obtain results on the behavior of boundary spectra under various operations on sheaves and, as a consequence, under integral transforms. We refer to section \ref{sec:sixoperations} for the precise statements. Here we only note that, as a corollary, one may obtain results on the local monodromy action on intersection cohomology (see Theorem \ref{thm:intersectioncomplex} and Corollary \ref{cor:intersectioncohom}.) 

\subsubsection{Monodromy in abelian covers}

As an application of our main Theorem \ref{thm:maintthm}, we obtain positive results towards Question \ref{question:questionintro} in some situations. For example, we show that the roots of unity can be obtained {\it uniformly} in $n$, and, moreover, our methods provide a schema to find such roots explicitly. We also obtain applications to the local monodromy of Alexander modules. We refer to Section \ref{sec:alexmodules} (in particular, Theorem \ref{thm:alexandermon}, Example \ref{example:alexmod}, and Corollary \ref{cor:abcovers}) for precise statements.

\subsection{Historical/Related work:}\label{sec:history/related}
We first note that Theorem \ref{thm:maintthm} (1) generalizes the classical monodromy theorem. Specifically, let $f: X \rightarrow S$ be as in the theorem (with $\dim(S) =1)$, $s \in S$, $\s{G}$ a constant local system and consider $\gamma \in {\rm AL}(S)$ centered at $s \in S$. In particular, $h: \pdisk \rightarrow S \setminus s$. Since $\spr(\gamma', \s{G}) = \{1\}$ for any $\gamma' \in {\rm AL}(X)$, the theorem shows that the local monodromy of $R^if_*(\s{G})$ is quasi-unipotent. In particular, this recovers (at least in char. 0), the classical local monodromy theorems of Brieskorn, Clemens, Grothendieck and Landman (\cite{SGA7-V1, Brieskorn, Clemens, Landman}). \\

In (\cite{SGA7-V1}), Grothendieck gives two proofs of the local monodromy theorem: one purely Galois theoretic and another based on the computation of his nearby cycles functor in the case where the special fiber is a divisor with normal crossings. The proofs of Theorems \ref{thm:maintthm} and \ref{thm:mainthm2} are a modification of the latter approach via nearby cycles. Analogous results in the context of variation of (mixed) Hodge structures (resp. regular singular connections) were obtained by Borel-Schmid (\cite{Schmid}) (resp. Katz (\cite{Katz})). We note that, in the Hodge theoretic setting, both Schmid and Katz obtain bounds on the level of quasi-unipotency in terms of the Hodge level. We do not obtain such bounds below. On the other hand, our results are applicable to sheaves of $R$-modules.\\

In \cite{Kashiwara-unip}, Kashiwara defined the notion of a quasi-unipotent constructible sheaf. More precisely, a constructible sheaf $\s{G}$ on $X$ is quasi-unipotent if $\spr(\gamma, \s{G})$ is contained in the set of roots of unity for all $\gamma \in {\rm AL}(X)$. In loc. cit., Kashiwara shows that (for proper morphisms) $Rf_*$ preserves the category of quasi-unipotent sheaves. We note that this is also a special case of our Theorem \ref{thm:maintthm}. Moreover, the results of this paper also prove the analogous assertion without any assumption on $f$ and also for $Rf_!$.\\

\subsection{Contents}
As noted above, the strategy for proving Theorems \ref{thm:maintthm} and \ref{thm:mainthm2} follows the strategy of Grothendieck's {\it geometric proof} of quasi-unipotence of local monodromy via nearby cycles. We briefly recall the contents of each section. \\

In Section \ref{sec:background}, we recall some basic background and set up some notation for the following sections. In Section \ref{sec:bmon}, we define various notions of {\it loops} and show that they give rise to the same spectra. In Section \ref{subsec:bmonoid} we recall the notion of boundary monoids and explain how Theorem \ref{thm:maintthm} (3) follows from \ref{thm:maintthm} (1). In Section \ref{sec:vancycles}, we recall some basic properties of nearby cycles. In Section \ref{sec:gpcohomology}, we discuss spectra in the setting of group actions on group cohomology. These will be crucially applied in Section \ref{sec:mainresults}.\\

In Section \ref{sec:mainresults}, we prove our main results when $\dim(S)=1$. In Section \ref{sec:mainresults1}, we use resolution of singularities arguments to reduce to a good setting, and give a basic vanishing cycles computation in the good setting. This proves the main theorem for $\dim(S) =1$, and also for $f_{!}$ in the case of the higher dimension. In Section \ref{sec:mainextensions}, we consider some natural extensions of Theorem \ref{thm:mainthm2} to a slightly more general setting, which will be useful in our applications to computing the monodromy of abelian covers.\\ 

In Section \ref{sec:dim>1}, we explain how to deduce Theorem \ref{thm:maintthm} (1) for $f_*$ from that of $f_!$ in the case where $\dim(S) \geq 1$. The crucial ingredient is to reduce to the setting of open immersions (in a good setting), and deal with this case directly (see Proposition \ref{prop:openimmersion}).\\

In sections \ref{sec:sixoperations} and \ref{sec:alexmodules} we give our applications to monodromy of integral transforms and monodromy of abelian covers, respectively.\\

{\bf Acknowledgements:} The authors thank the referee for a careful reading of the paper and for helpful comments and suggestions that improved it.\\

{\bf Notation:} In the following, $R$ will denote a commutative noetherian ring of finite global dimension and $K = \bar{K}$ will denote an algebraically closed field. For a complex algebraic variety $X$, $\D^b_c(X,R)$ denotes the bounded derived category of constructible sheaves of $R$-modules on $X$; if $R=K$, we denote this by $\D^b_c(X)$.\\

\section{Preliminaries}\label{sec:background}

\subsection{Remarks on Algebraic Monodromy}\label{sec:bmon}
In this section, we recall some equivalent characterizations of analytic loops and boundary spectra. Recall that if $X$ is a complex algebraic variety,  $\D^b_c(X)$ denotes the bounded derived category of $K$-vector spaces (where $K$ is an algebraically closed field).\footnote{Here constructible means in the underlying complex analytic topology but with stratifications given by  locally closed subsets.} Note that in the following, by abuse of notation, we will often view $X$ as a complex analytic space and simply use the same notation $X$ for the associated complex analytic space $X^{an}$.\\

We first explain how to associate an analytic map $h: \pdisk \rightarrow X$ to an analytic loop $\g \in {\rm AL}(X)$, and that given a constructible complex $\s{G}$ on $X$ one can choose $\pdisk$ small enough so that $h^*(\s{G})$ is locally constant.\footnote{A constructible complex $\s{G}$ is {\it locally constant} if all its homology sheaves are locally constant.} In order to see this, first note that the image of $\g$ is contained in an affine open $\Sp(A) \subset X$. Since $A$ is a finitely generated $\C$-algebra, we may choose a presentation $A = \C[x_1,\ldots,x_k]/(f_1,\ldots,f_s)$. A morphism $\g: \Sp(F^{an}) \rightarrow \Sp(A)$ is given by a collection of elements $H_1,\dots,H_k \in F^{an}$ which satisfy the polynomials $f_i$. Each $H_i$ defines a holomorphic function on a small punctured disk, and therefore we are given $k$ holomorphic functions on a small punctured disk. Since these satisfy the polynomials $f_i$, they give rise to an analytic map $h: \pdisk \rightarrow X$. Let $Z$ denote the  closure of $\gamma$ in $\Sp(A)$. Since $\s{G}|_Z$ is constructible, it is locally constant outside of a subvariety $W \subset Z$. Therefore, we may assume that $h^*{\s{G}}$ is locally constant after possibly shrinking the disk. Note that $\spr(\g,\s{G})$ is independent of the chosen disk (and the presentation). The discussion here also applies to constructible sheaves of $R$-modules, with ${\rm Sp}(\g,\s{G})$ defined as a closed subscheme of $\mathbb{A}^1_R$ (as given in Definition \ref{defn:spRmod}).

\bd\label{defn:curveloops}
\begin{enumerate}
\item Given $\s{G} \in \D^b_c(X)$ and $\gamma \in {\rm AL}(X)$, let $$\spr(\gamma, \s{G}) := \bigcup_{i \in \Z} \spr(\gamma,\s{H}^i(\s{G})).$$ Note that, since $\s{G}$ has bounded cohomology, only finitely many sets appear in the above union. 
Similarly, for $\s{G} \in \D^b_c(X,R)$, let 
$${\rm Sp}(\gamma, \s{G}) := \sum_{i \in \Z} {\rm Sp}(\gamma,\s{H}^i(\s{G})).$$
 
\item We set ${\rm BSp}(\s{G}) = \bigcup_{\gamma \in {\rm AL}(X)} \spr(\gamma, \s{G}). $
\end{enumerate}
\ed

We may also define {\it algebraic and formal} loops and consider analogously defined spectra. 

\bd
\begin{enumerate}
\item An element of $\gamma \in {\rm AL}(X)$ is an {\it algebraic loop} if the morphism 
$\gamma: \Sp(F^{an}) \rightarrow X$ factors through $\Sp(L)$ where $L$ is a field of transcendence degree at most one (over $\C$). 
\item Let $\C((t))$ denote the field of Laurent series. A {\it formal loop} $\gamma$ is a morphism
of schemes $\gamma: \Sp(\C((t))) \rightarrow X$ over $\C$.
\item Let $\s{G} \in \D_{c}^b(X)$. We set ${\rm BSp}_{a}(\s{G})$ to be the union of ${\rm Sp}_{red}(\gamma,\s{G})$ over all algebraic loops $\gamma$.
\item We say that a loop $\gamma \in {\rm AL}(X)$ is constant if it factors through a field of transcendence degree zero (i.e. $\Sp(\C)$).
\end{enumerate}
\ed

\begin{rmk}
Consider a pair $(T \subset \bar{T})$ where $T$ is a smooth curve, $\bar{T}$ is a smooth compactification of $T$, and a morphism $h: T \rightarrow X$. Then every point $s \in \bar{T} \setminus T$ gives rise to an algebraic loop of $X$. On the other hand, every algebraic loop arises in this manner.
\end{rmk}

\begin{rmk}
 Note that ${\rm BSp}_a(\s{G}) \subset {\rm BSp}(\s{G})$.
\end{rmk}

We may also define the eigenvalues of monodromy along formal loops as follows. Let ${\rm FL}(X)$ denote the set of formal loops, and $\gamma \in {\rm FL}(X)$. Let $Y$ denote the  closure of the image of $\gamma$, and $U \subset Y$ be the smooth locus of the local constancy locus of $\s{G}|_Y$.
Let $\bar{U}$ be a smooth compactification of $U$ with complement $D$ given by a divisor with simple normal crossings. By the valuative criterion of properness, we may now extend $\gamma$ to $\Sp(\C[[t]])$, that is, one has a commutative diagram (of $\C$-schemes):
$$
\xymatrix{
\Sp(\C((t))) \ar[r] \ar[d] &  U \ar[d] \\
\Sp(\C[[t]]) \ar[r] & \bar{U} }
$$
Let $x_0 \in \bar{U}$ be the image of the closed point $s_0 \in \Sp(\C[[t]])$; there exists a chart $(z_1,\ldots,z_n)$ around $x_0$ such that $D$ is given by $z_1 \cdots z_r =0$ for some $r \geq 0$. We may consider the pullback of $z_i$, and these can be written as $u_it^{k_i} \in \C[[t]]$ where $u_i$ is a unit. Retaining the $n$-tuple $(k_1,k_2,...,k_n)$ from the previous sentence, and given a tuple $\wt{u}:= (\wt{u_1},\ldots,\wt{u_n})$ of units of the ring of convergent power series, gives a holomorphic map $\delta$ from a small disk
$\Delta$ to a neighborhood of $x_0$ in $\bar{U}$ given by $z_i=\wt{u_i}t^{k_i}$
 for $1\leq i\leq n$. These $\delta$ (for varying choices of $\wt{u}$) restrict to holomorphic maps $\Delta^*\to U$, and moreover different choices of tuples $\wt{u}$ give maps which are homotopic to each other. We \emph{define} $\t{Sp}_{red}(\g)=\t{Sp}_{red}(\delta)$. 
Note that the analytic loop $\delta$ is algebraic if $\wt{u_i}=1$ for all $i$, and the $z_1,z_2,...,z_n$ belong to the co-ordinate ring of a  neighborhood of $x_0$ in $\bar{U}$. 
Moreover, the definition $\t{Sp}_{red}(\g)=\t{Sp}_{red}(\delta)$ is easily checked to be independent of the choice of the compactification $\bar{U}$ of $U$. We gather the results of this passage in the following lemma.

\bl\label{lem:allspsame}
With notation as above:
\begin{enumerate}
\item Given $\gamma \in {\rm FL}(X)$, $\spr(\gamma, \s{G})$ is independent of the choice of $U$ and compactification. We denote by 
${\rm BSp}_f(\s{G}) := \bigcup_{\gamma \in {\rm FL}(X)} \spr(\gamma, \s{G})$.
\item Given an analytic loop $\gamma$, there is an algebraic loop $\delta$ such that $${\rm Sp}_{red}(\gamma,\s{G}) = {\rm Sp}_{red}(\delta,\s{G}).$$
\item ${\rm BSp}_f(\s{G}) = {\rm BSp}(\s{G}) = {\rm BSp}_{a}(\s{G})$
\end{enumerate}
\el
\begin{proof}
We give details for part (2). The last part follows from part (2) and the remarks above. Let $\gamma$ be an analytic loop. We fix a compactification $\bar{U}$ as above. Let
$z_1,\dots,z_n$ be a system of parameters for the local ring
$\mathcal O_{\overline U,x_0}$ at $x_0$. Note that analytically the $z_i$ give local coordinates near $x_0$.
Let $u_1 t^{k_1},\dots,u_n t^{k_n}\in \mathbb C\{t\}$ be the images of the $z_i$ under the morphism $\mathcal O_{\overline U,x_0}\to \mathbb C\{t\}$ induced by $\gamma$, where each $u_i$ is a unit. In the chart given by the $z_i$ near $x_0$, the morphisms
\[
t\mapsto (u_1 t^{k_1},\dots,u_n t^{k_n})
\]
and
\[
\gamma' : t\mapsto (t^{k_1},\dots,t^{k_n})
\]
are homotopic. Hence $\gamma$ and $\gamma'$ have the same reduced spectrum. Now $\gamma'$ corresponds to a morphism of rings $\mathcal O^{an}_{\overline U,x_0}\to \mathbb C\{t\}$
sending $z_i$ to $t^{k_i}$. Restricting this morphism to $\mathcal O_{\overline U,x_0}\subset \mathcal O^{an}_{\overline U,x_0}$
gives a morphism $\mathcal O_{\overline U,x_0}\to \mathbb C\{t\}$ sending $z_i$ to $t^{k_i}$.
Since $\operatorname{Frac}(\mathcal O_{\overline U,x_0})$
has transcendence degree $n$ over $\mathbb C$, and the
$z_1,\dots,z_n$ form a transcendence basis,
every element of $\mathcal O_{\overline U,x_0}$
is algebraic over $\mathbb C(z_1,\dots,z_n)$.
Therefore, every element in the image is algebraic over
$\mathbb C(t)$.
Hence the image factors through a subfield $L\subset \mathbb C\{t\}$
which is algebraic over $\mathbb C(t)$.
Since $L$ has transcendence degree $1$ over $\mathbb C$,
the factorization $\mathcal O_{\overline U,x_0}\to L$
defines an algebraic loop $\delta$ having the same reduced spectrum
as $\gamma$. 
Finally, we leave the details of the first part to the reader.
\end{proof}

We record the following lemma for future use.

\bl\label{lem:basicproperties}
Let $X$ be a complex algebraic variety, and $\gamma \in {\rm AL}(X)$.
\begin{enumerate}
\item Given an exact triangle
 $$\s{F} \rightarrow \s{G} \rightarrow \s{H} \rightarrow \s{F}[1]$$
 in $\D^b_c(X,R)$, one has ${\rm Sp}(\gamma, \s{G}) \subset {\rm Sp}(\gamma,\s{F}) + {\rm Sp}(\gamma,\s{H})$, ${\rm Sp}(\gamma, \s{F}) \subset {\rm Sp}(\gamma, \s{G})$, and ${\rm Sp}(\gamma, \s{H}) \subset {\rm Sp}(\gamma, \s{G})$.
 \item Given an exact triangle
 $$\s{F} \rightarrow \s{G} \rightarrow \s{H} \rightarrow \s{F}[1]$$
 in $\D^b_c(X,K)$, one has $\spr(\gamma, \s{G}) \subset \spr(\gamma,\s{F}) \cup \spr(\gamma,\s{H})$. It follows that ${\rm BSp}(\s{G}) \subset {\rm BSp}(\s{F}) \cup {\rm BSp}(\s{H})$.The reverse inclusions hold for an exact sequence of sheaves.
 \item If $f: Y \rightarrow X$   is a morphism of complex algebraic varieties, and $\s{G} \in \D^b_c(X,R)$, then 
 ${\rm Sp}(\gamma', f^*\s{G}) = {\rm Sp}(f \circ \gamma' ,\s{G})$ where $\gamma' \in {\rm AL}(Y)$.
 \item If $f: Y \rightarrow X$   is a morphism of complex algebraic varieties, and $\s{G} \in \D^b_c(X)$, then 
 $\spr(\gamma', f^*\s{G}) = \spr(f \circ \gamma' ,\s{G})$ where $\gamma' \in {\rm AL}(Y)$. It follows that ${\rm BSp}(f^*\s{G}) \subset {\rm BSp}(\s{G})$. 
 
\end{enumerate}
\el
\bp
The proofs are standard and left to the reader.
\ep
Finally, we also record the following for future use:

\bl\label{lem:spdisjointunion}
Let $S = \coprod S_i$ be a finite disjoint union of schemes, and $\s{G} \in \D^b_c(S,R)$. Then:
\begin{enumerate}
\item ${\rm AL}(S) = \bigcup_i {\rm AL}(S_i)$,
\item For $\gamma \in {\rm AL}(S_i), {\rm Sp}(\gamma, \s{G}|_{S_i}) = {\rm Sp}(\gamma,\s{G})$, and 
\item if $R=K$, then ${\rm BSp}(\s{G}) = \bigcup_i {\rm BSp}(\s{G}|_{S_i}).$ 
\end{enumerate}
\el
\bp
The first two parts follow directly from the definition, and the last statement is a consequence of the definition of ${\rm BSp}$ and the first two parts.

\ep
\subsection{Boundary Monoid}\label{subsec:bmonoid}
Let $X$ be a smooth variety and $D = \bigcup_{i=1}^k D_i$ be a simple normal crossings divisor (s.n.c.d) where $D_i$ are the irreducible components. For a subset $I \subset \{1,\ldots,k\}$, $D_I := \bigcap_{i \in I} D_i$.\\

Given a point $x \in X$, there is a chart $U$ around $x$ where, if $z_1,\ldots,z_n$ are the local coordinates, then $U \cap D$ is given by $z_1\cdots z_d =0$. In this case, the fundamental group
$\pi_1(U \setminus D)$ is the free abelian group generated by the canonical loops (i.e. the counterclockwise loops resulting from identifying $U \setminus D$ with $(\Delta^{*})^{d}\times \Delta^{n-d}$) around each of the $D_i$. The image of the submonoid generated by these loops in $\pi_1(X \setminus D)$ is independent of the choice of $U$ (and the base points) up to conjugacy. In particular, it gives rise to a commutative monoid in $\pi_1(X \setminus D)$, well defined up to conjugacy. Moreover, up to conjugacy, a different choice of $x$ in the same irreducible component of $ D_I \setminus \bigcup_J D_J$ (where the union is over $J$ such that $|J| > |I|$) gives the same monoid up to conjugacy. In particular, one has a finite number of commutative monoids well defined up to conjugacy in $\pi_1(X \setminus D)$. \\

Let $Y$ be a smooth complex algebraic variety, $\s{L}$ be a local system of $K$-vector spaces on $Y$, and $\bar{Y}$ a smooth compactification of $Y$ with $D: = \bar{Y} \setminus Y$  an s.n.c.d. If $\gamma \in {\rm AL}(Y)$, and $h: \pdisk \rightarrow Y$ is the associated analytic map, then $h$ extends to a complex analytic map $\bar{h}: \disk \rightarrow \bar{Y}$. If $\bar{h}(0) \in Y$, then $\spr(\gamma, \s{L}) = 1 \in K^{\times}$. Otherwise, let $\bar{h}(0) = x \in D$ denote the center of the disk. In this case, one obtains an element of the monoid given by the local fundamental group of $x$ defined above. As a consequence, we note that ${\rm BSp}(\s{L}) = \bigcup_{(i,\gamma \in M_i)} \spr(\gamma, \s{L})$ (where $M_i$ are the finite number of commutative monoids obtained as in the previous paragraph). In particular, if $U$ is a curve, then 

$${\rm BSp}(\s{L}) = \{1\} \cup (\bigcup_i \bigcup_{n \in \mathbb{N}} \spr(\gamma_i^n,\s{L}))$$ 

where $\gamma_i$ are the loops around the boundary points. 

\begin{rmk}\label{rmk:bmonoid}
Note that the discussion in the previous paragraph shows that the image of ${\rm AL}(Y)$ to the set of the conjugacy classes of $\pi_1(Y)$,
is the image of the union of the boundary monoids in the same set. Note that the set of boundary monoids depends on the chosen compactification $\bar{Y}$ of $Y$.
\end{rmk}

The discussion above allows us now to give a direct proof of Theorem \ref{thm:maintthm} (3) {\it assuming} Theorem \ref{thm:maintthm} (1).

\begin{proof}(Theorem \ref{thm:maintthm} (1) implies Theorem \ref{thm:maintthm} (3))
Let $f: X \rightarrow S$ and $\s{G}$ be as in the theorem. First, note that we may assume that $X$, $S$ are reduced and connected. Consider $\s{H}:=R^qf_{?}(\s{G})$. Let $S = \coprod S_i$ where each $S_i$ is locally closed smooth and $\s{L}_i: = \s{H}|_{S_i}$ is a local system. Let $S_i \subset \bar{S}_i$ be smooth compactifications with boundary a simple normal crossings divisor. By the above discussion, we have, for each $i$, a finite number of boundary monoids $M_i^{j}$ such that
$${\rm BSp}(\s{L}_i) = \bigcup_{(j,\gamma \in M_i^{j})} {\rm Sp}_{red}(\gamma,\s{L}_i) = \bigcup_{(j,\gamma \in M_i^{j})} {\rm Sp}_{red}(\gamma,\s{H}).$$ On the other hand (by \ref{lem:spdisjointunion}) $${\rm BSp}(\s{H}) = \bigcup_i {\rm BSp}(\s{L}_i).$$

Note that each monoid is finitely generated and commutative. Since the monoid is commutative, we see that the eigenvalues of the action of a loop in a particular monoid are given by products of eigenvalues of the generators of that monoid. Since there are only finitely many generators and finitely many monoids, we may apply \ref{thm:maintthm} to each of these finitely many loops. In particular, we take for $r'$ the lcm of the integers $r$ obtained for each such loop via Theorem \ref{thm:maintthm} (1).
\end{proof}

\subsection{Nearby Cycles}\label{sec:vancycles}

We recall some standard results on the nearby cycles functor. We refer to (\cite{KS}, Section 8.6) or (\cite{SGA7-V2}) for more details. Let $X$ be a complex algebraic variety, $f: X \rightarrow S$ a morphism to a smooth curve, and $\disk$ a small disk centered at a point $s_0 \in S$. By abuse of notation, we denote by $f:X \rightarrow \disk$ the restriction of $f$ to the disk $\disk$. Let $\s{G} \in \D^b_c(X;R)$, and $\pdisk$ the disk with the origin (i.e. $s_0$) removed. Let $\pi: \widetilde{\pdisk} \rightarrow \pdisk$ denote the universal cover. Explicitly, we consider the map $p: \Delta \rightarrow \Delta^*$ with $z \mapsto e^{2\pi i z}$.\footnote{We shall ignore base points in the discussion below.} Now consider the resulting Cartesian diagram:

$$
\xymatrix{
X_0 \ar[r]^{i} & X  \ar[d] & \widetilde{X} \ar[d] \ar[l]^{\tilde{j}}\\
  & \disk   &  \widetilde{\pdisk} \ar[l] }
$$

Here $X_0:= f^{-1}(0)$, and the nearby cycles complex is defined as follows:

$$R\Psi_{f}(\s{G}):= i^*R\tilde{j}_*\tilde{j}^*\s{G}.$$ 

The natural deck transformation $T: \widetilde{\pdisk} \rightarrow \widetilde{\pdisk}$ (corresponding to the canonical generator of $\pdisk$) gives rise to the {\it monodromy morphism}

$$T: R\Psi_f(\s{G}) \rightarrow R\Psi_f(\s{G}).$$

We recall some basic properties of the nearby cycle functor:
\begin{enumerate}
\item $R\Psi_f(\s{G})$ is a constructible complex on $X_0$.
\item If $f$ is proper (and up to further shrinking of the disk), then there is a spectral sequence with 

$${\rm E}_{2}^{p,q} := \H^p(X_0,R^q\Psi_f(\s{G})) \Rightarrow \H^{p+q}(X_t,\s{G})$$

where $t \in \pdisk$ is a general point 
\item Note that for $\disk$ small, $\H^{i}(X_t,\s{G})$  is a local system on $\pdisk$, and in particular comes equipped with a monodromy action. The aforementioned spectral sequence is compatible with the monodromy actions.
\item Let $x \in X_0$. The stalk $R^i\Psi_f(\s{G})_x$ can be computed as follows. Let $B(x,\varepsilon)$ be an open ball of radius $\varepsilon$ centered at $x$ in $X$. Then for all $0 < \varepsilon << 1$ and $0< \delta << \varepsilon$ the aforementioned stalk can be identified with $\H^i(B(x,\varepsilon) \cap f^{-1}(t), \s{G})$ where $0< |t| < \delta$.

\item 
Let $T: \s{F} \rightarrow \s{F}$ be a morphism of sheaves of $R$-modules on $X$. In this case, we may view $\s{F}$ as a sheaf of $R[x]$-modules (or $R[x,x^{-1}]$-modules, with $x$ acting via $T$), and consider its annihilator $Ann(\s{F}) \subset R[x]$. Let ${\rm Sp}(T,\s{F}) \subset \A^1_R$ denote the corresponding closed subscheme. We may also work point-wise and define the annihilators $Ann(\s{F}_y)$ for $y \in X$.  Applying $R\Gamma(X,-)$ gives rise to a functor with values in the bounded derived category of $R[x]$-modules. In this setting, one has:
$$Ann(\oplus R^i\Gamma(X,\s{F})) \supset Ann(\s{F}) \subset Ann(\s{F}_y). $$
Note that if $\s{F}$ is a constructible sheaf of $R$-modules, then $R\Gamma(X,\s{F})$ is an object of the bounded derived category of finitely generated $R$-modules. 
    
\end{enumerate}

\subsection{Group Cohomology}\label{sec:gpcohomology}
Let $G$ be a group, $H \subset G$ a normal subgroup, and $G':= G/H$. Let $M$ be a $G$-module, where $M$ is a (finitely generated) $R$-module and the $G$-action is $R$-linear. \\

Given $g \in G$, one has an induced $R$-linear map $\rho_g: M \rightarrow M$, and as before one may view $M$ as an $R[x]$-module and consider the corresponding scheme theoretic support ${\rm Sp}(g,M) \subset \Sp(R[x])$, and the corresponding ideal $I(g,M) \subset R[x]$.\\

With $G$ and $H$ as above, the exact sequence 

$$1 \rightarrow H \rightarrow G \rightarrow G' \rightarrow 1$$

gives rise to an ($R$-linear) action of $G'$ on the (group) cohomology groups $\H^i(H,M)$. We briefly recall a description of this action and refer to (\cite{Brown}, III.8) for the details. Given $g \in G$, let $c_g: H \rightarrow g^{-1}Hg = H$ denote the map $c_g(h) = g^{-1}hg$. The maps $c_g$ and $\rho_g: M \rightarrow M$ induce a map on cohomology $\H^i(H,M) \rightarrow \H^i(H,M)$ as follows. We may view the domain of $\rho_g$ as an $H$-module where $h \in H$ acts via $\rho \circ c_g(h) = \rho(g^{-1}hg)$; with this modified action on the domain, $\rho_g$ is a morphism of $H$-modules. An application of the usual bi-functoriality of group cohomology (contravariant in the first variable and covariant in the second) gives the desired morphism (of $R$-modules) $\H^i(H,M) \rightarrow \H^{i}(H,M)$. If $g \in H$, then this map is trivial, and therefore the action factors through $G'$. 

\bl\label{lem:gpcohomsp}
Let $H,G,G',R$ be as above, and let $M$ be a finitely generated $R$-module with an $R$-linear $G$-action. Let $g' \in G'$, and $g \in G$ be a lift of $g'$. Suppose $g \in Z(G)$. Then $I(g,M) \subset I(g',\H^i(H,M))$.
\el
\bp
Let $p(x) \in I(g,M)$. It follows that $p(g)$ annihilates the $R$-module $M$. We would like to show that $p(g')$ annihilates $\H^i(H,M)$. Let $Z$ be the center of $R[G]$. It is enough to show that if $\rho(z)$ annihilates $M$, then the action of $\rho(z')$ (with $z'$ the image of $z$ in $R[G']$) on $\H^i(H,M)$ is also trivial.\\

We may regard $\H^i(H,-)$ as a functor from $R[G]$-modules to $R[G']$-modules. 
Now $\rho(z)$ is a $G$-module endomorphism of every $R[G]$-module $M$, and therefore induces the natural transformation $\H^i(H,\rho(z))$
from the functor $\H^i(H,-)$ to itself. As above, let $z'$ denote the image of $z$ in $R[G']$; we have the action of  $\rho'(z')$ on $\H^i(H,M)$. Their difference $\H^i(H,\rho(z))-\rho'(z')$
is a natural transformation from the functor $ff\H^i(H,-)$
to itself, where $ff$ denotes the forgetful functor from $R[G']$-modules to $R$-modules. This difference is zero on $ff\H^0(H,-)$.  The system of $ff\H^i(H,-)$ forms a sequence of effaceable cohomological $\delta$-functors, so their difference is zero for all $i\geq 0$.
In particular, if $M$ is an $R[G]$-module for which $\rho(z)=0$, then the action of $\rho'(z')$ on $\H^i(H,M)$ is also trivial.

\ep

\begin{rmk}\label{rem:gpcohomspfield}
If $R = K$ is an algebraically closed field, $V$ is a finite dimensional $K$-vector space and $G$ is abelian, then the previous lemma shows that the eigenvalues of $g' \in G'$ acting on $\H^i(H,V)$ are contained in the eigenvalues of a lift $g \in G$ of $g' \in G'$ acting on $V$.
\end{rmk}

\section{Proof of Theorem \ref{thm:maintthm} (1), (2) and Theorem \ref{thm:mainthm2}: $\dim(S)=1$.} \label{sec:mainresults}
In this section, we prove our main results when $\dim(S) =1$. In the first subsection, we reduce the statements to a {\it good situation} (see below) and give an explicit computation of the monodromy action of stalks of nearby cycles in the good situation. The main theorems are deduced from this result. In the second subsection, we prove some generalizations to the setting where $R$ is replaced by a locally constant sheaf of $R$-modules. The latter result will be useful in the application to monodromy of abelian covers in Section 6.

\subsection{Reductions and Key Proposition}\label{sec:mainresults1}
In this section, we prove our main results in the setting where $\dim(S) =1$. More precisely, we prove Theorem \ref{thm:maintthm} (1), (2) and Theorem \ref{thm:mainthm2} in the setting where $\dim(S) =1$. We shall deduce all three statements by reducing it to the following setting.\\

Let $f: X \rightarrow S$ be a morphism with $\dim(S) =1$, and $\s{G} \in \D^b_c(X)$ (or $\D^b_c(X,R)$). We say that $(X,S,f,\s{G})$ is in the {\it good situation} if the following holds:
\begin{enumerate}
\item We have a commutative diagram:
$$
\xymatrix{
U \ar[r]^{j} \ar[dr] & X  \ar[r]^{\bar{j}}  \ar[d]^{f} & \bar{X} \ar[d]^{\bar{f}} \\
   & S  \ar[r] & \bar{S}  }
$$
where all the horizontal arrows are open immersions, $\bar{f}$ is proper, all the varieties in the diagram are smooth (connected), $\bar{S},\bar{X}$ are proper, and $D: = \bar{X} \setminus U$ is an s.n.c.d. \item Moreover, $D = A + B$ where $A,B$ are s.n.c.d.'s with no common components, $\bar{X} \setminus X =A$, and $X \setminus U = B \setminus A \cap B$. 
\item For $s \in \bar{S}$, let $\bar{X}_s$ denote the corresponding scheme theoretic fiber. Then $(\bar{X}_s)_{red} \cup D$ is an s.n.c.d. for all $s \in \bar{S}$. We have a local system $\s{L}$ on $U$, and $\s{G}:= j_!\s{L}$. 

\end{enumerate}

\bt\label{thm:reductions}
Suppose that Theorems \ref{thm:maintthm} (1), (2), and \ref{thm:mainthm2} hold for all quadruples $(X,S,f,\s{G})$ in the good setting. Then, Theorems \ref{thm:maintthm} (1), (2) (resp. \ref{thm:mainthm2}) hold for all quadruples $(X,S,f,\s{G})$ where $f: X \rightarrow S$ is a morphism with $\dim(S)=1$ and $\s{G} \in \D^b_c(X)$ (resp. $\s{G} \in \D^b_c(X,R)$).
\et
\begin{proof}
Let $f: X \rightarrow S$ be a morphism of 
schemes of finite type over $\C$, with $\dim(S) =1$ and $\s{G} \in \D^b_c(X)$. We note that all the reductions below are also valid in the setting of Theorem \ref{thm:mainthm2} and $\s{G} \in \D^b_c(X,R)$. We begin with some preliminary reductions:
\begin{enumerate}
\item[1:] First, note that we may assume that all the schemes in question are connected and reduced. 
\item[2:] We may assume that the morphism $f$ is dominant. Otherwise, the image is a collection of points, and the local monodromy on a zero dimensional scheme is trivial.
\item[3:] We may assume that the base $S$ is a smooth connected curve. Let $\tilde{S} \rightarrow S$, denote the normalization. To see this, note that since $\dim(S) =1$, the natural map
$${\rm AL}(\tilde{S}) \rightarrow {\rm AL}(S)$$ is a bijection for non-constant loops since morphisms from $\Sp(F^{an}) \rightarrow S$ factor through the generic point. 
\end{enumerate}

We are now in the setting where $f: X \rightarrow S$ is a morphism of complex algebraic varieties with $S$ a smooth connected curve, and $f$ is dominant. We shall proceed via induction on $\dim(X)$. 

{\bf Step 1:} Since $\s{G}$ is constructible,  there is an open dense subset $j: U \hookrightarrow X$ such that $\s{L} := \s{G}|_U$ is locally constant. Moreover, up to replacing $U$ by a smaller  open subset, we may assume $U$ is smooth. Let $Z = X \setminus U \xhookrightarrow{i} X$ denote the closed complement. Then one has the standard triangle
$$
j_!\s{L} \rightarrow \s{G} \rightarrow i_*i^*\s{G}
$$
in $\D^b_c(X)$. By applying induction and Lemma \ref{lem:basicproperties} (since the complement $Z: = X \setminus U$ is of lower dimension), it is sufficient to prove the theorem for $j_!\s{L}$. In particular, we now assume that $j: U \hookrightarrow X$ is an open immersion with $U$ smooth, $\s{L}$ is a local system on $U$, and $\s{G} = j_!\s{L}$. Note that we may replace $U$ by an open dense $V \subset U$. This will be necessary in Step 4 below.\\
{\bf Step 2:} Consider a commutative diagram

$$
\xymatrix{
U \ar[r]^{\tilde{j}} \ar[dr]^{j} & Y \ar[r]^{\tilde{f}} \ar[d]^{\pi} & S \\
    &     X \ar[ur]^{f}  & 
}
$$

where $j$, $\tilde{j}$ are open immersions and $\pi$ is a proper morphism. The diagram above induces (via functoriality of sheaves on $X$-schemes) the isomorphism: 

$$R\pi_{!}\tilde{j}_!(\s{L}) \rightarrow j_!(\s{L})$$

Moreover, since $\pi$ is proper, one has $R\pi_! = R\pi_*$. Because $R\widetilde{f}_?=Rf_? \circ R\pi_?$ we obtain the isomorphism below for 
$?=!$ and $?=*$:

$$R\tilde{f}_{?}(\tilde{j}_!(\s{L})) \rightarrow Rf_{?}(j_!\s{L}).$$

Note that $\pi$ induces a map ${\rm AL}(Y) \rightarrow {\rm AL}(X)$ so that, if $\gamma' \in {\rm AL}(Y)$ and $U \times_X Y \rightarrow U$ is an isomorphism, then 
${\rm Sp}(\gamma',\tilde{j}_!\s{L}) ={\rm Sp}(\pi \circ \gamma',j_!\s{L}) $. It follows that it's enough to prove the theorem for the morphism $\tilde{f}$ and $\s{G} = \tilde{j}_!\s{L}$. \\
{\bf Step 3:} We apply the previous step to an embedded resolution of singularities of the pair $(X,Z)$. In particular, we may assume that $X$ is smooth and $Z$ is an s.n.c.d. Moreover, we may choose a smooth compactification $\bar{X}$ of $X$ such that $\bar{X} \setminus U = D$ is an s.n.c.d. We may further assume (once again applying resolution of singularities) that $D = A + B$ where $A$ and $B$ are s.n.c.d.'s and $\bar{X} \setminus X = A$. In particular, $X \setminus U = B  \setminus A \cap B$. \\
{\bf Step 4:} We fix a smooth compactification $S \subset \bar{S}$ so that the morphism $f$ extends to $\bar{f}: \bar{X} \rightarrow \bar{S}$. For $s \in \bar{S}$, let $\bar{X}_s$ denote the corresponding scheme theoretic fiber. Applying relative desingularization to a certain pair $(\bar{X},D')$ over $\bar{S}$, we may assume that for every $s \in S$, $(\bar{X}_s)_{red} \cup D$ is an s.n.c.d. Specifically, suppose $(D_i)_{i \in I}$ are the irreducible components of $D$. For each $J \subset I$, we restrict $\bar{f}$ to $D_J = \bigcap_{j \in J} D_j$, and consider the restriction $\bar{f}_J$ of $\bar{f}$ to $D_J$. We apply Sard's theorem to the morphisms $\bar{f}_J$ to obtain an  open subset $S' \subset S$ over which each $\bar{f}_J$ is smooth. We now consider $Z:= \bar{X} \setminus \bar{f}^{-1}(S')$, set $D' := Z \cup D$, and apply resolutions to the pair $(\bar{X},D')$ to obtain a pair $(\bar{\bar{X}},\bar{D})$ where $\bar{D}$ is total transform of $D'$. Note that our original $U$ is now replaced by $U_{S'} := U \times_S S'$, and $X$ by its inverse image in $\bar{\bar{X}}$. Finally, note that new $(X,S,f,\s{G})$ is now in the good situation.
\end{proof}

We shall now assume that $(X,S,f,\s{G})$ is in the good situation and we are concerned with the boundary monodromy of $R^if_{?}\s{G}$. We first consider the case of $? = *$. In fact, the case of $? =!$ is simpler and will follow immediately from the method of proof of the former case. \\

Consider now $\gamma \in {\rm AL}(S)$, and the associated morphism from the punctured disk $h: \Delta^{*} \rightarrow S $. Since $\bar{S}$ is compact, this extends to a map from the disk $h: \Delta \rightarrow \bar{S}$. Seeing $\gamma$ as a loop in ${\rm AL}(\bar{S})$, one has 

$${\rm Sp}(\gamma,R^if_*\s{G}) = {\rm Sp}(\gamma, R^i\bar{f}_*(R\bar{j}_*\s{G}))$$
if $\s{G}$ is a constructible sheaf of $R$-modules and because $R\bar{f}_*R\bar{j}_*(\s{G})|_S = Rf_*(\s{G})$. The analogous claim also holds in the setting of sheaves of $K$-vector spaces. 
Let $\s{H} := R\bar{j}_*\s{G}$. The well-known results listed in section \ref{sec:vancycles} show that ${\rm Sp}_{red}(\gamma, R^i\bar{f}_*(\s{H}))$ is controlled by the monodromy action on the stalks of nearby cycles of $\s{H}$. The latter action is dealt with in the proposition that follows.  \\

The above assumptions have as a consequence the following normal form for the diagram above (restricted to $\Delta$). Let $x_0 \in \bar{f}^{-1}(0)$. In this case, one has a chart $\Omega$ around $x_0$ (and centered at $x_0$) with local coordinates given by 

$$a_1,\ldots,a_{\ell},a'_{1},\ldots,a'_{\ell'},b_1,\ldots,b_m,b'_1,\ldots, b'_{m'},c_1,\ldots,c_n,c'_1,\ldots,c'_{n'},$$ 

and $N := \ell+\ell'+n+n'+m+m'$ such that 
\begin{enumerate}
\item $A \cap \Omega$ (resp. $B \cap \Omega$) is the divisor defined by the vanishing of $\prod_{i=1}^{\ell} a_i \prod_{j=1}^{\ell'} a'_{j}$ (resp. $\prod_{i=1}^{m} b_i \prod_{j=1}^{m'} b'_{j}$).
\item The morphism $\bar{f}$ is given by $\bar{f}(a_1,\ldots) = a_1^{\lambda_1}...a_{\ell}^{\lambda_{\ell}}b_1^{\mu_1}...b_m^{\mu_m}c_1^{\nu_1}...c_n^{\nu_n}$ for 
positive integers $\lambda_i, \mu_j,\nu_k$.
\end{enumerate}
The neighborhood $\Omega$ can be identified with a product of small disks $\Delta^N$ (with coordinates as above), and with this notation $X \cap \Omega$ is a product 

$$(\Delta^{*})^{\ell + \ell'} \times \Delta^{N- (\ell + \ell')}$$

and 
$U$ is the product

$$(\Delta^{*})^{\ell + \ell' + m + m'} \times \Delta^{N'}$$

where $N' = N - (\ell + \ell' + m + m')$. 
With this notation, let $p: \Delta^N \rightarrow \Delta^{N''}$ denote the projection to the {\it non-primed} coordinates (so that $N'' = \ell+m+n$) and $g: \Delta^{N''} \rightarrow \Delta$ the map given by $g(a_1,\ldots,a_{\ell}, b_1,\ldots,b_m,c_1,\ldots,c_n) = a_1^{\lambda_1}...a_{\ell}^{\lambda_\ell}b_1^{\mu_1}...b_m^{\mu_m}c_1^{\nu_1}...c_n^{\nu_n}$. In particular, $\bar{f} = g \circ p$ (when restricted to $\Omega$). We are interested in the eigenvalues of the monodromy action on the stalk $R\Psi_{\bar{f}}(\s{H})_{x_0}$ of the nearby cycles along the morphism $\bar{f}$ (restricted to a small disk $\Delta$ via $h$) of the sheaf $\s{H}$. 

\begin{prop}\label{prop:mainprop}
With notation as above, and suppose $R = K$ an algebraically closed field:
\begin{enumerate}
\item Suppose that $x_0$ is contained in an irreducible component of $B$ which is not contained in $\bar{f}^{-1}(0)$ (i.e. $m' > 0$), then $R\Psi_{\bar{f}}(\s{H})_{x_0} = 0$
\item Suppose $m' = 0$, and $n > 0$. Let $\gamma \in \pi_1(\Delta^{*})$ be the canonical generator as before. The action of $\gamma$ on $R^i\Psi_{\bar{f}}(\s{H})_{x_0}$ has roots of unity as eigenvalues. Moreover, there is a $\tilde{\g} \in {\rm AL}(U)$ and an integer $ r > 0$ such that $f \circ \tilde{\g} = \g^r$, the eigenvalues of
$\g$ acting on $R\Psi_{\bar{f}}(\s{H})_{x_0}$ are contained in ${\rm Sp}_{red}(\tilde{\g},\s{L})^{\frac{1}{r}}$, and the action of $\tilde{\g}$ is trivial. Finally, the action of $\g$ is diagonalizable if the characteristic of $K$ is zero.
\item Suppose $m' = 0$ and $n = 0$. In this case, there is either an $a$ or a $b$ variable. Then there is a $\tilde{\g} \in {\rm AL}(U)$ and an integer $r>0$ such that $f \circ \tilde{\g} = \g^r$ and the eigenvalues of
$\g$ acting on $R\Psi_{\bar{f}}(\s{H})_{x_0}$ are contained in ${\rm Sp}_{red}(\tilde{\g},\s{L})^{\frac{1}{r}}$.
\end{enumerate}
\end{prop}
\bp
Before beginning the proof, we set-up some notation. By the discussion in \ref{sec:vancycles}, the stalk $R\Psi_{\bar{f}}(\s{\s{H}})_{x_0}$ can be computed as follows.
For $t \in \Delta^{*}$ small enough, the stalk is given by (and choosing $\Omega$ small enough) the cohomology group $\H^i(\Omega \cap \bar{f}^{-1}(t), \s{H}) = \H^i(\Omega \cap f^{-1}(t), \s{G}) $. In the proof below, by abuse of notation, we shall still use $U, X, \bar{X}$ to denote $U \cap \Omega, X \cap \Omega$ and $ \bar{X} \cap \Omega = \Omega$ (and similarly for the fiber $f^{-1}(t)$).
\begin{enumerate}
\item With notation as above, $f^{-1}(t)=g^{-1}(t)\times (\Delta^{*})^{\ell'} \times \Delta^{m'} \times  \Delta^{n'} $ where $(\Delta^{*})^{\ell'}$ (resp. $\Delta^{n'}$, $\Delta^{m'}$) is the product of punctured disks (resp. disks) in the $a'$ variables (resp. $b',c'$ variables).
Now note that $f^{-1}(t)\cap B$ is the same as above, except that the product of the disks in the $b'$ variables is replaced by its closed subset $b'_1b'_2...b'_{m'}=0$. On the other hand, $\s{G}$ restricted to $B$ is zero (since it is by definition $j_!\s{L}$). It now follows by Corollary \ref{cor:homotopycor} that the stalk of the nearby cycles vanishes when $m'>0$. We apply the corollary as follows. We take $W=\Omega\cap f^{-1}(t)$, we express $W=W_1\times W_2$ where $W_2$ is the product of the discs in the
$b'_j$ variables, and define $F:W\times I\to W$ by $F(w_1,b'_1,b'_2,...,b'_{m'},s)=(w_1, sb'_1,sb'_2,...,sb'_{m'})$. 

\item Suppose $m' =0$. Let $G$ denote the fundamental group of the open subset given by 
$a_1...a_lb_1...b_mc_1...c_na'_1...a'_{l'}\neq 0$. This set is the product of punctured disks in the $a$,$a'$,$b$,$c$ variables and the disks in the $c'$ variables. Let $V$ denote this set, and we shall view it as an open subset of $U$ (or rather $U \cap \Omega$). Note that $U$ is the product of punctured disks in the $a$,$a'$,$b$, variables and the disks in the $c$ and $c'$ variables. Note that one has open inclusions:

$$ V \subset U \subset X \subset \bar{X}.$$

Now $G$ is the free abelian group with basis given by the canonical loops in the punctured disks. We denote these by 
\[\g_i(a),\g_j(b),\g_k(c),\g_{i'}(a')\text{  where }1\leq i\leq l,1\leq j\leq m,1\leq k\leq n,1\leq i'\leq l'.\]
The inclusion of $V \subset U$ gives a surjection
$$G \to G'',$$
where $G''$ is $\pi_1(U)$, and the map is given by killing the loops in the $c$ variables. On $U$, the constructible sheaf $\s{G}$ is the local system $L$ and therefore corresponds to a finite dimensional representation, also denoted by $L$, of $G''$. As noted above, we are interested in the monodromy action on $\H^{i}(f^{-1}(t),\s{G})$. 
First, note that $t \neq 0$, and in particular $f^{-1}(t) \subset U$ since $m'=0$ and therefore all $\mu_i > 0$ (for all $1 \leq i \leq m$) or there are no $b$-variables (in which case $U = X$). In particular, $\s{G}$ is the local system given by the local system $L$ on $f^{-1}(t)$. Similarly, since none of the $c$-variables can vanish on a point in $f^{-1}(t)$, it follows that $f^{-1}(t) \subset V$. Now, the morphism $\bar{f}$ restricted to $V$ induces a morphism 
$\eta: G \rightarrow \g^{\Z}$ on fundamental groups, and we let $K=\ker(\eta)$. We are interested in the action of $\g$ on $\H^i(f^{-1}(t), L).$ Note that $K$ is the fundamental group of a connected component of $f^{-1}(t)$. The number of such connected
components is $d$ where $\eta(G)=\g^{d\Z}$. Note that $\eta(\gamma_i(a)) = \g^{\lambda_i}$, $\eta(\gamma_j(b)) = \g^{\mu_j}$, $\eta(\gamma_k(c))= \g^{\nu_k}$ and $\eta(\gamma_{i'}(a')) = 0$. It follows that $d$ is given by the g.c.d. of the $\lambda_1,\ldots,\lambda_l$, $\mu_1,\ldots,\mu_m$, and $\nu_1,\ldots,\nu_n$. Below, let $r_{\lambda_i} = \lambda_i/d$, and define $r_{\mu_j}, r_{\nu_k}$ similarly. We shall consider two cases: $d = 1$, and $d$ arbitrary.
\begin{enumerate}
\item[Case (i)] Suppose $d =1$. In this case, $\eta$ is surjective, and $\H^i(f^{-1}(t),L)$ can be identified with the group cohomology $\H^i(K,L)$ (since $f^{-1}(t)$ is an Eilenberg-MacLane space with fundamental group $K$). Moreover, the monodromy action can be identified with the action of the quotient $G/K$ on the aforementioned cohomology group. It now follows from Lemma \ref{lem:gpcohomsp} and Remark \ref{rem:gpcohomspfield} that 
the eigenvalues of $\gamma$ can be computed by choosing a lift of $\gamma$ to $G$. If $n > 0$, then there is a $c$ variable appearing; we may assume that $\nu_1 >0$. In particular, the action of $\gamma^{\nu_1}$ is given by the induced action of $\gamma_1(c)$. But, the latter acts trivially, and therefore $\gamma^{\nu_1}$ acts as the identity. It follows that the eigenvalues of $\gamma$ are contained in the $\nu_1$-th roots of unity. In this case, the $\g$-action is diagonalisable if the characteristic of $K$ does not divide $\nu_1$.
\item[Case (ii):] Suppose $d >1$. Again, since there is a $c$ variable, we may assume $\nu_1 > 0$. In this case, we have the action of $G/K\cong\g^{d\Z}$ on $\H^i(K,L)$ and therefore the natural action of $\g^{\Z}$ on $\t{Ind}^{\g^{\Z}}_{\g^{d\Z}}\H^i(K,L)$. Moreover, this representation (of $\g^{\Z}$) can be identified with the monodromy action of $\g$ on $\H^i(f^{-1}(t), \s{G})$. One can see this geometrically as follows. Consider the morphism $f' = f^{1/d}$ (which is well defined in the current setting), and consider the resulting cartesian diagram (in the neighborhood $\Omega$ of $x_0$):
$$
\xymatrix{
\Omega \ar[r]^{Id} \ar[d]^{\bar{f'}} & \Omega  \ar[d]^{\bar{f}} \\
\pdisk \ar[r]^{z \mapsto z^d} &  \pdisk .}
$$
Restricting to $X$ (and denoting by $f'$ the resulting morphism), and arguing as in the previous case, we see that the assertion holds for $f'$ and the sheaf $\s{G}$. In particular, $\gamma'^{\nu_1/d}$ is the identity on $\H^i(K,L)$, where $\gamma'$ denotes the loop in the disk on the left in the diagram above; it maps to $\gamma^{\nu_1}$ in the disk on the right. Note that $R^if_*(\s{G}) \cong p_*R^if'_*(\s{G})$ where $p(z)= z^d$, and (by generic base change) the stalks, for small $t$, of these sheaves are the cohomology of the fibers. It follows that the local system $\H^i(f^{-1}(t),\s{G})$ (on the disk on the right) is identified 
with the push-forward of the local system $\H^{i}(f'^{-1}(t), \s{G})$ (on the disk on the left) along the morphism $z \mapsto z^d$. In particular, it is identified with the induced representation above. 
The action of $\g^d$ on $\H^i(f^{-1}(t), \s{G})$ is  therefore the natural action of $\g^d$ on the direct sum of $d$ copies
of $\H^i(K,L)$. It follows that $\gamma$ has roots of unity as eigenvalues. More precisely, $(\gamma)^{\nu_1}$ is the identity.

\end{enumerate}
\item Suppose now that both $m'=0$ and $n=0$. We may once again argue as in the second case above. In all cases (including those considered in (2) above), we may consider $f' = f^{1/d}$ with $d$ as above.
Note that $\eta(\g_i(a))=\g^{\lambda_i}$. Suppose that there is an $i$ such that $\lambda_i>0$ (i.e. there is an $a$-variable). Arguing as in (2), the eigenvalues of $\g$ are therefore contained in the set of $\lambda_i$-th roots of the eigenvalues
of $\g_i(a)$ on $V$. If there are no $a$-variables, then there is certainly a $b$-variable, and we see that the eigenvalues of $\g$ on the nearby
cycles are contained  in the $\mu_1$-th roots of the eigenvalues of $\g_1(b)$. 
\end{enumerate}
We note that in all cases considered in (2) the eigenvalues of $\gamma$ are roots of unity, and the $\g$-action is diagonalizable if the characteristic of $K$ is 0 (or more generally if the characteristic does not divide any of the $\lambda_i, \mu_j, \nu_k$). 
\ep

\bl\label{lem:homotopylemma}
Let $W$ be a paracompact Hausdorff space, $I = [0,1]$, and $\s{H}$ be a sheaf of abelian groups on $W \times I$ such that 
\begin{enumerate}
\item $\s{H}_{{w} \times (0,1]}$ is locally constant for all $w \in W$,
\item $\s{H}|_{W \times \{0\}} = 0$.
\end{enumerate}
Then $\H^{i}(W \times I, \s{H}) = 0$ for all $i \geq 0$.
\el
\bp
Suppose $W$ is a single point. Then the statement amounts to the assertion that if $\s{H}$ is a sheaf on $I$ which is locally constant sheaf on $(0,1]$ and whose stalk at $0$ is $0$, then $\s{H}$ is cohomologically acyclic. Therefore, the claim is clear in this case. 
Let $p: W \times I \rightarrow W$ denote the projection map. By proper base change and previous discussion, $Rp_*(\s{H}) = 0$. Consideration of the Leray spectral sequence now gives the desired conclusion. 
\ep

\bcor \label{cor:homotopycor}
Let $W$ and $I$ be as in Lemma \ref{lem:homotopylemma}. Let $F: W \times I \rightarrow W$ be a continuous map such that $F(w,1) = w$ for all $w \in W$, $\s{G}$ be a sheaf on $W$, and $\s{H} = F^{*}(\s{G})$. Suppose $\s{H}$ satisfies the hypotheses of Lemma \ref{lem:homotopylemma}. Then $\H^{i}(W,\s{G}) = 0$ for all $i \geq 0$.
\ecor
\bp
Let $i: W \rightarrow W \times I$ denote the map $i(w) = (w,1)$. Then $F \circ i = Id_W$, and therefore the induced composite morphism
$$
\H^{i}(W, \s{G}) \xrightarrow{F^*} \H^i(W \times I, \s{H} ) \xrightarrow{i^*} \H^i(W,\s{G})
$$
is the identity map. On the other hand, by Lemma \ref{lem:homotopylemma}, the middle term is zero. Therefore, $\H^{i}(W, \s{G}) =0$ for all $i \geq 0$.

\ep

\bp(Theorem \ref{thm:maintthm} (1), $\dim(S) =1$)
By Theorem \ref{thm:reductions}, it is enough to prove this in the good situation. We explain how to deduce this from Proposition \ref{prop:mainprop}.
\begin{enumerate}
\item[(i)] We first consider the case of $Rf_*$. By the discussion preceding Proposition \ref{prop:mainprop}, it is enough to prove that there is a finite set $M$ of loops on $X$  so that $${\rm Sp}_{red}(\gamma,R^i\bar{f}_*\s{H})\subset \bigcup_{\gamma' \in M} {\rm Sp}_{red}(\gamma', \s{G}).$$ In fact, we shall see that there is a finite set $M \subset {\rm AL}(U)$ with the requisite property.
Recall that $R\Psi_{\bar{f}}(\s{H})$ is a constructible complex of $K[T,T^{-1}]$-modules. Fix a stratification of $\bar{f}^{-1}(0)$ on which $R\Psi^i_{\bar{f}}(\s{H})$ is locally constant, and let $Z$ be a connected component of one of the strata. By Proposition \ref{prop:mainprop}, given a point $x_0$ in this stratum, either the stalk $R^i\Psi_{\bar{f}}(\s{H})_{x_0} = 0$, or there is a loop $\tilde{\gamma} \in {\rm AL}(U)$ such that 
$f \circ \tilde{\gamma} = \gamma^r$ and $${\rm Sp}_{red}(\gamma,R^i\Psi_{\bar{f}}(\s{H})_x) \subset {\rm Sp}_{red}(\tilde{\gamma},\s{L})^{\frac{1}{r}}.$$ It follows that 
$${\rm Sp}_{red}(\gamma, R^i\Psi_{\bar{f}}(\s{H})|_Z) \subset {\rm Sp}_{red}(\tilde{\gamma},\s{L})^{\frac{1}{r}}.$$
Since there are only a finite number of strata, we conclude that there is a finite subset $M \subset {\rm AL}(U)$ such that for all $\gamma' \in M$, $f \circ \gamma' = \gamma^r$, and 
$$
{\rm Sp}_{red}(\gamma, R^i\Psi_{\bar{f}}\s{H}) \subset \bigcup_{\g' \in M} {\rm Sp}_{red}(\g',\s{L})^{\frac{1}{r}}.
$$
Note that we can choose an $r$ and a finite set $M$ that works for all $i$. By the discussion in Section \ref{sec:vancycles} (5), we have 
$${\rm Sp}_{red}(\gamma, \H^j(\bar{f}^{-1}(0),R^i\Psi_{\bar{f}}(\s{H})) \subset {\rm Sp}_{red}(\gamma, R^i\Psi_{\bar{f}}\s{H}).$$
We now apply the nearby cycles spectral sequence (see section \ref{sec:vancycles} (2)). Note that the abutment is $$\H^p(\bar{X}_t,\s{H}) = R^p\bar{f}_*(\s{H})_t.$$ 
Since this is an extension of subquotients of $\H^j(\bar{f}^{-1}(0),R^i\Psi_{\bar{f}}\s{H})$, the result follows.

\item[(ii)] Consider now the case of $Rf_{!}(\s{G})$. Note that in this case, we are reduced to computing ${\rm Sp}_{red}(\g, R\bar{f}_*(\bar{j}_!\s{G})$. Therefore, this is the good situation where $X = \bar{X}$ and $A = \emptyset$.
\end{enumerate}

\ep

\bp(Theorem \ref{thm:maintthm} (2))
This is a consequence of the reductions above. The key point is that on an open subset of $S$, $R^if_*\s{G}$ (resp. $R^if_!\s{G})$) is a local system. Since $\dim(S) =1$, the complement is a finite set. In particular, there are only a finite number of algebraic loops to consider on the base $S$ (the ones with center in this finite set and those with center in $\bar{S} \setminus S$). More precisely, the discussion in Section \ref{subsec:bmonoid} and Remark \ref{rmk:bmonoid} combined with the first part of the theorem now immediately gives the second part. 

\ep

\bp(Theorem \ref{thm:mainthm2}, $\dim(S) =1$)
Suppose now that $\s{G}$ is a constructible sheaf of $R$-modules. Again, by Theorem \ref{thm:reductions}, we may assume that we are in the good situation. Arguing as in the proof of Theorem \ref{thm:maintthm} (1) above, it is enough to prove the analog of Proposition \ref{prop:mainprop} in the setting of sheaves of $R$-modules. Note that in the last step of the argument in the proof of Theorem \ref{thm:maintthm} (1) (for $Rf_*$) when passing to extensions, the union of reduced spectra must be replaced by the {\it sum} of the corresponding subschemes. We now explain the modifications needed for the analog of Proposition \ref{prop:mainprop} in the setting of $R$-modules. The proof of part (1) of that proposition clearly goes through in the setting of $R$ modules, giving the same conclusion in the case that $m' >0$. One can deal with parts (2) and (3) of the proposition simultaneously as follows. First, we set $f' = f^{1/d}$ as in the proof of Proposition \ref{prop:mainprop}, and continue with the notation above. Now suppose $P(T^{r_{\lambda_i}})$ (resp. $P(T^{r_{\mu_j}})$, $P(T^{r_{\nu_k}})$) is a polynomial that annihilates $\H^i(f'^{-1}(t),V)$, then the induced representation described in the proof of Proposition \ref{prop:mainprop} is annihilated by $P(T^{\lambda_i})$ (resp. $P(T^{\mu_j})$, $P(T^{\nu_k}$)). This gives the desired result.

\ep

\subsection{Some remarks and extensions} \label{sec:mainextensions}
In this section, we discuss some extensions of Theorem \ref{thm:mainthm2} to a slightly more general setting. These will be useful in section \ref{sec:alexmodules} for our applications to monodromy in abelian covers.\\

Let $f: X \rightarrow S$ be a morphism with $\dim(S) =1$ and consider the following data:
\begin{enumerate}
\item A {\it locally constant} sheaf $\s{R}_S$ of commutative noetherian rings (of finite homological dimension) on $S$.
\item A sheaf $\s{F}$ of $f^{-1}\s{R}_S$-modules on $X$ which is {\it weakly constructible} as a sheaf of abelian groups. Recall, this means that there is a good stratification on which $\s{F}$ is locally constant, but we do not require any {\it finiteness} hypotheses.
\item A locally constant sheaf of ideals $\s{I}_S \subset \s{R}_S$,
\end{enumerate}

We first note that the functors $Rf_*$, $Rf_!$ are still defined on such objects. Let $s_0 \in S$, and consider a loop $\gamma \in {\rm AL}(S)$ with center $s_0$. In particular, the morphism $h: \Delta^* \rightarrow S$ associated to $\gamma$ extends to the full disk $h: \Delta \rightarrow S$ so that $h(0) = s_0.$ Upon restricting everything to this disk we have the following data:
\begin{enumerate}
\item $\s{R}$ (resp. $\s{I}$) is constant (up to shrinking $\Delta$), and canonically identified with its stalk $R := \s{R}_{s_0}$ (resp. $I:= \s{I}_{s_0}$).  
\item $R^if_{?}(\s{F})$ is a weakly constructible sheaf of $R$-modules, and up to shrinking the disk is locally constant on $\pdisk$. Note that since $\s{R}$ is locally constant on the disk, the monodromy action on $R^if_{?}(\s{F})$ is $R$-linear. We are interested in ${\rm Sp}(\gamma, R^if_{?}(\s{F}))$.
\item Let $\s{R}_{\s{I}} : = \s{R}/\s{I}$. Note that this is also a locally constant sheaf of commutative noetherian rings (of finite homological dimension). Again, we may assume that its restriction to the disk is constant and canonically identified with the stalk $R/I$. Let $\s{F}_{\s{I}} := \s{F} \otimes_{f^{-1}(\s{R})} f^{-1}(\s{R}/\s{I})$. 
\item $R^if_{?}(\s{F}_{\s{I}})$ is a weakly constructible sheaf of $R/I$-modules, and we are also interested in ${\rm Sp}(\gamma, R^if_{?}(\s{F}_{\s{I}}))$.
\end{enumerate}

We claim that the conclusion of Theorem \ref{thm:mainthm2} remains valid in the previous setting. In addition, we also have a compatibility property when going modulo $I$. Given a closed subscheme $Z \subset \A^1_R$ and $I \subset R$, we denote by $Z_I := Z \bigcap \A^1_{R/I}$ the corresponding scheme theoretic intersection considered as a closed subscheme of $\A^1_{R/I}$.

\begin{thm} \label{thm:mainthmvariant}
With notation as above, there is a finite set $M$ of pairs $(\g',n_{\g'})$ where $\g' \in {\rm AL}(X)$ such that:
\begin{enumerate}
\item We have 
$${\rm Sp}(\g,R^if_{?}(\s{F}) \subset \sum_{\gamma' \in M} {\rm Sp}(\g',\s{F})^{[1/n_{\g'}]}.$$
\item We have
$${\rm Sp}(\g,R^if_{?}(\s{F}_I) \subset \sum_{\gamma' \in M} {\rm Sp}(\g',\s{F})^{[1/n_{\g'}]}_I.$$
\end{enumerate}

\end{thm}

\begin{proof}
This follows from the following observations, whose details we leave to the reader:
\begin{enumerate}
\item Firstly, for both assertions, the reductions of the previous section to the good Hironaka situation can be performed in our given setting.
\item Secondly, once in the good situation, we note that locally around the loop, we are dealing with constructible sheaves of $R$-modules. In particular, this is exactly the setting of (the proof of) Proposition \ref{prop:mainprop}. This immediately proves the first assertion. 
\item We assume that we are in the local setting of Proposition \ref{prop:mainprop}. Our sheaf $\s{F}$ is the $\s{G}$ of loc. cit. and $\s{F}_I$ is the sheaf $\s{G}_I$. Consider a stratification such that both $R^i\Psi_f(\s{F})$ and $R^i\Psi_f(\s{F}_I)$ are locally constant. Given a stratum $Z$ (of such a stratification), we see that for an $x \in Z$:
$${\rm Sp}(\g, R^i\Psi_f(\s{F})_x) \subset {\rm Sp}(\g', \s{F})$$ and therefore
$${\rm Sp}(\g, R^i\Psi_f(\s{F})|_Z) \subset {\rm Sp}(\g', \s{F}).$$
We also have analogs of these inclusions for the sheaf $\s{F}_{I}$.
The stalks $R^i\Psi_f(\s{F})_x$ can be computed as in the proof of Proposition \ref{prop:mainprop}. We see that in each case it is either zero, or a certain group cohomology with coefficients in the module $M$ or $M/I$. Since $Ann_{R/I}(M) = (I + Ann_R(M))/I$, the result follows. 

\end{enumerate}    
\end{proof}

\section{Proof of Theorem \ref{thm:maintthm} (1) and Theorem \ref{thm:mainthm2}: $\dim(S) \geq 1$.} \label{sec:dim>1}

In this section, we will complete the proof of Theorems \ref{thm:maintthm} (1) and \ref{thm:mainthm2}. Let $f: X \rightarrow S$, and $\s{G} \in {\rm D}^b_c(X,R)$. \\

First note that, in order to prove Theorems \ref{thm:maintthm} (1) and \ref{thm:mainthm2} for $R^qf_{!}(\s{G})$, we may reduce to the case of $\dim(S) =1$ since any algebraic loop lies in a curve and since $Rf_!$ commutes with base change (see Lemma \ref{lem:allspsame}). Therefore, the claims for $Rf_!$ follow from the results of the previous section.\\

We consider the case for $Rf_*$ in a special setting. Suppose $\bar{X}$ is smooth and proper, and let $U \hookrightarrow X \hookrightarrow \bar{X}$ be open immersions where $U$ is  dense. Let $D := A + B \subset \bar{X}$ be a simple normal crossings divisor such that $\bar{X} \setminus A = X$, $\bar{X} \setminus D = U$, $U = X \setminus B \cap X$ and $A,B$ have no common components. We shall refer to such a triple $(U,X,\bar{X})$ as a `good Hironaka triple'.

\begin{prop}\label{prop:openimmersion}
Let $U \hookrightarrow X$ be an open dense subset and $X \hookrightarrow \bar{X}$ such that $(U,X,\bar{X})$ is a good Hironaka triple.
Let $L$ be a local system of $R$-modules on $U$ and $\gamma \in {\rm AL}(\bar{X})$. Let $j: U \hookrightarrow X$ and $j': X \hookrightarrow \bar{X}$ denote the given open immersions, with $A,B,$ and $D$ as above and consider $R^qj'_*(j_!L)$. Then there is a loop $\gamma' \in {\rm AL}(U)$ such that
${\rm Sp}(\gamma,R^qj'_*(j_!L)) \subset {\rm Sp}(\gamma', L) $.
If $L$ is a local system of $K$-vector spaces, it follows that ${\rm BSp}(R^qj'_*(j_!L)) \subset {\rm BSp}(L)$.
\end{prop}

Before proving the proposition, we prove Theorems \ref{thm:maintthm} (1) and \ref{thm:mainthm2} for $Rf_*$ assuming the proposition.

\begin{proof}(\ref{thm:maintthm} (1) and \ref{thm:mainthm2} for $Rf_*$)
We begin with some reductions. Recall, Theorem \ref{thm:mainthm2} implies Theorem \ref{thm:maintthm} (1) by passing to the underlying reduced scheme. Therefore, we shall only consider the former setting. Let $(X,S,f,\s{G})$ be as in Theorem \ref{thm:mainthm2}. Without loss of generality, we may assume that $\s{G}$ is a constructible sheaf (rather than a bounded complex of such) of $R$-modules and let $\gamma \in {\rm AL}(S)$. \\

{\bf Step 0:} As before, we may assume that $X$ and $S$ are connected and reduced. We may also assume that $S$ is proper.\\
{\bf Step 1:} We note that if Theorem \ref{thm:mainthm2} holds for morphism $g: Y \rightarrow Z$ and $h: Z \rightarrow S$, then it also holds for $h \circ g$ as an application of the Leray spectral sequence. We may factor our given morphism $f: X \rightarrow S$ as $X \xhookrightarrow{j} \bar{X} \xrightarrow{\bar{f}} S$ where the first morphism is an open immersion and the second is proper. Since $R\bar{f}_{!} = R\bar{f}_*$ for proper morphisms, we are reduced to proving \ref{thm:mainthm2} of an open immersion.\\
{\bf Step 2:} It remains to prove Theorem \ref{thm:mainthm2} for the open immersion $j: X \hookrightarrow \bar{X}$. We may stratify $X$ by smooth locally closed subsets $S_i \xhookrightarrow{j_i} X$ such that the restriction of $\s{G}$ to each of these is a local system. As before, this gives rise to a filtration of $\s{G}$ with associated graded of the form $j_{i,!}\s{G}|_{S_i}$. As an application of Lemma \ref{lem:basicproperties}, we are reduced to proving the claim for each such stratum. In particular, we may assume that there is a connected smooth locally closed subset $Z \subset X$, and $\s{G}$ is obtained as an extension by zero from a local system $L$ on $Z$.\\
{\bf Step 3:} Let $X'$ (resp. $\bar{X}'$) denote the closure of $Z$ in $X$ (resp. $\bar{X}$). Consider the following diagram:
$$
\begin{tikzcd}[row sep=1.6em, column sep=1.6em]
Z \ar[r, hook, "j_1"] \ar[dr, hook, "i"'] 
  & X' \ar[r, hook, "j'"] \ar[d, hook, "a"'] 
  & \bar{X}' \ar[d, hook, "b"] \\
  & X \ar[r, hook, "j"] 
  & \bar{X}
\end{tikzcd}
$$
By construction, $j_1$, $j$, $j'$ are open immersions, both $a,b$ are closed immersions, and the square is cartesian. It is enough to prove the Theorem for the open immersion $j'$ and the sheaf $\s{G}' := j_{1,!}L$.
To see this, note that ${\rm Sp}(\gamma, Rj_*(\s{G})) = {\rm Sp}(\gamma, Rj_*((a \circ j_1)_{!}L)) = {\rm Sp}(\gamma, R(j\circ a)_*(\s{G}')) = {\rm Sp}(\gamma, R(b\circ j')_*(\s{G}'))$. Since $b$ is a closed immersion, the latter is trivial if the image of $\gamma$ 
is in the complement of $\bar{X}'$. In particular, we may assume that $\gamma$ lands in $\bar{X}'$; in that case, we have 
${\rm Sp}(\gamma, R(b\circ j')_*(\s{G}')) = {\rm Sp}(\gamma, Rj'_*(\s{G}'))$. In particular, we are reduced to the setting where $U \subset X$ is an smooth open dense subset, $L$ is a local system on $X$, and $\s{G}$ is extension by zero of the local system $L$ on $U$.\\
{\bf Step 4:} We may further assume that $\bar{X}$ is proper by an application of Nagata compactification.  \\
{\bf Step 5:} An application of embedded resolutions allows one to obtain a commutative diagram:
$$
\xymatrix{
U \ar[d]^{=} \ar[r]^{j_U'} & V \ar[r]^{j'} \ar[d]^{\pi} & \bar{V} \ar[d]^{\bar{\pi}} \\
U \ar[r]^{j_U} & X \ar[r]^{j} & \bar{X} 
}
$$
where $(U,V,\bar{V})$ is a good Hironaka triple, the right square is cartesian, and the vertical maps are proper. Let $\s{G}' := j'_{U,!}L$, and note that $R\bar{\pi}_*(Rj'_*(\s{G}')) = R(j \circ \pi)_*(\s{G}') = Rj_*(\s{G})$. Therefore, Theorem \ref{thm:mainthm2} for the open immersion $j$ and sheaf $\s{G}$ is a consequence of Theorem \ref{thm:mainthm2} for the proper morphism $\bar{\pi}$ and Proposition \ref{prop:openimmersion}. In order to obtain a diagram as above, one proceeds as follows. First, let $D_0 := \bar{X} \setminus U$, and apply embedded resolution to the pair $(\bar{X},D_0)$ to obtain a proper birational map $\bar{\pi}_1: \bar{V}_1 \rightarrow \bar{X}$ (isomorphism over $U$) so that $\bar{V}_1$ is smooth proper and $D_1 : = \bar{\pi}_1^{-1}(D_0)_{red}$ is an sncd. Let $V_1$ be the pull back of $X$, and set $A_1 := \bar{V}_1 \setminus V_1$. We once again apply embedded resolutions to the pair $(\bar{V}_1,A_1)$(with boundary $D_1$) in order to obtain a proper birational morphism $\bar{\pi}_2: \bar{V} \rightarrow \bar{V}_1$ (isomorphism over $U$) so that $A := \bar{\pi}_2^{-1}(A_1)_{red}$ is an sncd, $D := \pi_2^{-1}(D_1)_{red}$ is an sncd, and $D = A + B$. We set $V$ to be the pull back of $V_1$, and note that the resulting $(U,V,\bar{V})$ is a good Hironaka triple. Finally, setting $\bar{\pi} := \bar{\pi}_1 \circ \bar{\pi}_2$ gives the desired diagram.
\end{proof}

In the remainder of this section, we give the proof of Proposition \ref{prop:openimmersion}. In particular, we now fix a good Hironaka triple $(U,X,\bar{X})$ with $D=A+B$ as above. If $\gamma \in  {\rm AL}(\bar{X})$, then the associated morphism $h: \Delta^* \rightarrow \bar{X}$ extends to the disk $h: \Delta \rightarrow \bar{X}$ since $\bar{X}$ is proper. By abuse of notation, we set $\gamma(0) := \bar{x} \in \bar{X}$ to be the image $h(0)$, and call this the center of $\gamma$. We begin with a lemma which considers the local analog of our setting and that will be useful in the proof.

\begin{lem}\label{lem:observation1}
Let $\Delta^n$ denote the open polydisk with coordinates $(z_1,\ldots,z_b,w_1,\ldots,w_m)$ 
with $n = b+m$, $b \geq 1$, $B := \{z_1 = \cdots = z_b =0\}$, $A := \{w_1 = \cdots = w_a = 0\}$ (for some $0 \leq a \leq m$), $X := \Delta^n \setminus A$, $U := \Delta^n \setminus (A + B)$, and $L$ a local system on $U$. Let $j: U \hookrightarrow X$, and $j': X \hookrightarrow \Delta^n$ denote the natural inclusions. Then $\H^q(\Delta^n,Rj'_*j_!L) = 0$ for all $q\geq0$.

\end{lem}
\begin{proof}
Let $z = (z_1,\cdots,z_b)$, $w = (w_1,\cdots,w_m)$, and $F: X \times [0,1] \rightarrow X$ denote the homotopy $(z,w,t)\mapsto(tz,w)$ for $(z,w)\in X$ and $t\in[0,1]$. Then $F$ satisfies:
\begin{enumerate}
\item[(a)] $F(u,1)=u$ for all $u\in X$
\item[(b)] $F(u,0)\in B \cap X$ for all $u\in X$
\item[(c)] The restriction of $F^*(j_!L)$ to $\{u\}\times (0,1]$ is locally constant for every $u\in X$
\end{enumerate}
Applying Corollary \ref{cor:homotopycor} to   $W=X$ and the sheaf $j_!L$, we deduce that $\H^q(\Delta^n,Rj'_*j_!L) = \H^q(X, j_!L)= 0$ for every $q\geq 0$. 
\end{proof}

\bp (Proof of Proposition \ref{prop:openimmersion})
We begin the proof with setting up some notation and making two preliminary observations. For $a\geq 0$, the set of $x\in\bar{X}\setminus B$ that belong to exactly $a$ irreducible components of $A$ is denoted by $T^a$. Note that $T^0 = U$. We define $T^{-1}=B \cap X$. Thus, $\bar{X}$ is the disjoint union of its  locally closed subsets $T^a$ taken over $a=-1,0,1,\ldots$. The proposition is concerned with the sheaves $R^qj'_*\s{G}$ where $\s{G}:=j_!L.$\\

We first observe that the sheaves $R^qj'_*\s{G}$ vanish for all $q\geq 0$ when restricted to $B.$ This is a consequence of Lemma \ref{lem:observation1}. By definition the stalk
$R^qj'_*\s{G}_x$ (for $x \in \bar{X}$) is the direct limit of $\H^q(\Omega\cap X,j_!L)$ taken over all neighborhoods $\Omega$ of $x$ in $\bar{X}$. 
So, it suffices to check that
$\H^q(\Omega\cap X,j_!L)=0$ for a cofinal system of neighborhoods $\Omega$ of $x$ in $\bar{X}$. On the other hand, if $x$ belongs to exactly $b\geq 1$ irreducible components of $B$, then $x$ has a fundamental system of neighborhoods of the form given in Lemma \ref{lem:observation1}.\\

We next note that by (\cite{Schurm}, 4.2.1) the sheaves $R^qj'_*\s{G}$ are \emph{locally constant} when restricted to $T^a$ for every $a \geq 0$.\footnote{The result in loc. cit. is in fact more general and with weaker assumptions.}  We give a self-contained proof here for the convenience of the reader and which makes explicit the local structure of this sheaf. It is a simple consequence of the product structure (in the usual topology)
induced by the stratification. More precisely, every point $x$ of $T^a$ (with $a \geq 0$) has 
\begin{enumerate}
\item[(i)] a neighborhood $\Omega \cong \Delta^a \times \Delta^{n-a}$ in $\bar{X}\setminus B$ such that $\Omega\cap X= (\Delta^*)^a \times \Delta^{n-a}$
\item[(ii)] a locally constant sheaf $L'$ on
$(\Delta^*)^a$ such that the pull-back of $L'$ to $\Omega \cap X$ along the projection map ($(\Delta^*)^a \times \Delta^{n-a} \rightarrow (\Delta^*)^{a}$) is isomorphic to the restriction of $L$ to $\Omega\setminus A$.
\end{enumerate}

By (ii), we have $$\H^q((\Delta^*)^a,L')\ \cong \H^q(\Omega\cap X,\s{G}).$$
From this we easily deduce that every $x'\in\Omega\cap T_a= 0 \times \Delta^{n-a}$
possesses a fundamental system of neighborhoods $\Omega'$ of $x'$ such that
$ \H^q(\Omega\cap X,\s{G})\to\H^q(\Omega'\cap X,\s{G})  $ is an isomorphism. This produces a natural isomorphism from the constant sheaf 
\begin{equation}\label{eqn:stalkcompute}
\H^q((\Delta^*)^a,L')_{\Omega\cap T_a} \xrightarrow{\cong} \t{R}^qj'_*\s{G}|_{\Omega\cap T_a}.
\end{equation}

The proposition will now be deduced from Lemma \ref{lem:gpcohomsp}. By Lemma \ref{lem:spdisjointunion}, an analytic loop $\g$ of $\bar{X}$, by virtue of being an $F^{an}$-valued point $\bar{X}$, is in fact an $F^{an}$-valued point of $T^a$ for
a unique $a$. The action of $\g$ on the local system 
$R^qj'_*\s{G}|_{T^a}$ is under discussion. By Lemma \ref{lem:observation1} and the discussion above, we only need to discuss the case $a\geq 0$. One the other hand, the case $a = 0$ is clear. From now on we assume that $a > 0$.\\

Note that the center $\g(0)$ of the loop lies in the intersection of exactly $a + s$ components of $A$ for some $s \geq 0$. Assume that the center $\g(0)$ lies in exactly $b$ components of $B$.
Now we have a neighborhood $\Omega$ of $\g(0)$ and a pointed isomorphism of $(\Omega, \g(0))$ with $(\Delta^n = \Delta^a \times \Delta^s \times \Delta^b \times \Delta^{n-a-s-b}, 0)$ such that (under this identification):
\begin{enumerate}
\item[(i)]$(\Delta^*)^{a+s} \times \Delta^b \times \Delta^{n-a-s-b} =\Omega\cap X$ 
\item $ (\Delta^*)^{a+s} \times (\Delta^*)^b \times \Delta^{n-a-s-b}=\Omega\cap U$ 
\item[(iii)] there is a local system $L'$ on 
$(\Delta^*)^a\times(\Delta^*)^s\times (\Delta^*)^b$ which pulls back to $L$ on $\Omega\cap U$. 
\end{enumerate}
Denote by $G_a,G_s, G_b$ the fundamental groups of 
$(\Delta^*)^a,(\Delta^*)^s,(\Delta^*)^b$ respectively. The fundamental group of 
$\Omega\cap U$ is given by the fundamental group of $(\Delta^*)^a\times(\Delta^*)^s\times(\Delta^*)^b$. Choose a point $t\in \Omega\cap U$ and let $M$ denote the stalk of $L'$ at the image of $t$. Thus, $M$ is an
$R[G]$-module where
$G:=G_a\times G_s\times G_b$. Let $W=\{0\}\times (\Delta^*)^s\times (\Delta^*)^b \subset \Delta^{a+s+b}$ and note that the $W \times \Delta^{n-a-s-b}=\Omega\cap T^a$.\\

Now, the locally constant sheaf we are mainly concerned with, namely
$R^qj'_*\s{G}|_{T^a\cap\Omega}$, is clearly the pullback of a 
 sheaf $\s{F}$ on $W$. The explicit description (\ref{eqn:stalkcompute}) of this sheaf given above shows that the stalk of $\s{F}$ at the image of $t$  is identified with $\H^q(G_a\times 1\times 1,M)$. The fundamental group of $\Omega\cap T^a$ is identified with the quotient group $Q:=G/(G_a\times 1 \times 1) \cong G_s \times G_b$.\\

Now $\g$ has its image $q\in Q$. We choose a lift $q'\in G$ of $q$, and then define $\g'$ to be the preimage 
of $q'$ under the isomorphism $\pi_1(\Omega\cap U)\to G$. All the groups in question are commutative, and so the proposition follows from an application of Lemma \ref{lem:gpcohomsp}.
\ep

\section{Integral Transforms and Intersection cohomology}\label{sec:sixoperations}
In this section, we collect some general results on the behavior of local monodromy
under various functors and give an application to the behavior of local monodromy under integral transforms and intersection cohomology.

\subsection{Local Monodromy under integral transforms}
Let $f: X \rightarrow Y$ be a morphism of schemes, and as before $\s{G} \in \D^b_c(X)$ a constructible complex of $K$-vector spaces. 

\bt
With notation as above, let $\s{H} \in \D^b_c(Y)$.
\begin{enumerate} 
\item One has ${\rm BSp}(f^*\s{H}) \subset {\rm BSp}(\s{H}).$ 
\item There is an integer $r >0 $ (depending on $f, \s{G}$) such that 
${\rm BSp}(Rf_*\s{G})^+ \subset {\rm (BSp}(\s{G})^+)^{\frac{1}{r}}.$ The similar assertion also holds for $Rf_!$.
\item Given $\s{F},\s{G} \in \D^b_c(X)$, 
${\rm BSp}(\s{F} \otimes \s{G}) \subset {\rm BSp}(\s{F}){\rm BSp}(\s{G})$. Here, the right-hand side is the set consisting of the products of elements in each of the sets.
\item Given $\s{F},\s{G} \in \D^b_c(X)$, 
${\rm BSp}(\s{H}om(\s{F},\s{G})) \subset {\rm BSp}(\s{F})^{-1}{\rm BSp}(\s{G})$. Here ${\rm BSp}(\s{F})^{-1}$ is the set of $\lambda$ such that $\lambda^{-1} \in {\rm BSp}(\s{F})$.
\end{enumerate}
\et
\bp
\begin{enumerate}
\item We have already taken note of the case of $f^*$ (see \ref{lem:basicproperties}).
\item This is the main result of the previous sections.
\item This follows from the standard fact that the eigenvalues of a tensor product of matrices consist of the products of the eigenvalues of each of the matrices. 
\item This follows from the fact that for finite dimensional vector spaces $\Hom(V,W) \cong V^* \otimes W$ (as representations of some $G$). Moreover, the eigenvalues for the dual representation are given by the inverses of the eigenvalues of the original representation.
\end{enumerate}
\ep

As an application, we compute the monodromy for various integral transforms.

\begin{cor}
Consider a diagram of schemes 
$$
\xymatrix{
  & Z \ar[dl]^p \ar[dr]^q &  \\
X &  &  Y}
$$
Let $\s{K} \in \D^{b}_c(Z)$, and consider the functor $I: \D^b_c(X) \rightarrow \D^b_c(Y)$ where $I(\s{G}) = q_*(p^*(\s{G}) \otimes \s{K}).$
Then there is an integer $r >0$ such that 
$${\rm BSp}(I(\s{G}))^{+} \subset (({\rm BSp}(\s{G}){\rm BSp}(\s{K}))^+)^{\frac{1}{r}}.$$
\end{cor}

For example, this applies to the usual Radon transform (or more generally Brylinski-Radon transform). In the case of the Radon transform, $X = \P^n$ is projective space and $Y = \check{\P}^n$ is the dual projective space. If $H \subset \P^n \times \check{\P}^n$ denotes the usual incidence correspondence, then $\s{K} := i_*K$ i.e., the direct image of the constant sheaf $K$ under the inclusion $i: H \hookrightarrow \P^n \times \check{\P}^n$. With this notation, the Radon transform
$\s{R}: \D^b_c(\P^n) \rightarrow \D^b_c(\check{\P}^n)$ is (up to shifts) by definition
$q_*(p^*(\s{G}) \otimes \s{K}).$
In this case, ${\rm BSp}(\s{K}) = {\rm BSp}(K) = \{1\}$.

\begin{cor}
With notation as above, there exits $r > 0$ such that
$${\rm BSp}(\s{R}(\s{G}))^+ \subset ({\rm BSp}(\s{G})^{+})^{\frac{1}{r}}.$$
In particular, $\s{R}$ preserves the full subcategory of quasi-unipotent sheaves.
    
\end{cor}

\subsection{Intermediate Extensions}
In this section, we discuss the monodromy of intermediate extensions of perverse sheaves and, in particular, intersection cohomology. We denote by $\s{P}(X)$ the category of perverse sheaves on $X$ (with coefficients in $R$ where $R$ is also assumed to be an artinian ring). \\

Given a locally closed immersion $j: U \hookrightarrow X$, one has the {\it intermediate extension}
$$ j_{!*}: \s{P}(U) \hookrightarrow \s{P}(X).$$
In this setting, we have the following result for spectra of intermediate extensions.

\begin{thm}\label{thm:intersectioncomplex}
Let $\gamma \in {\rm AL}(X)$ and $\s{G} \in \s{P}(U)$. Then there is a finite set $M$ of pairs $(\gamma', n_{\g'})$ with $\gamma' \in {\rm AL}(U)$, $n_{\gamma'}$ a positive integer and such that 
$${\rm Sp}(\gamma,j_{!*}(\s{G})) \subset  \sum_{(\gamma',n_{\gamma'}) \in M} {\rm Sp}(\gamma',\s{G})^{[1/n_{\gamma'}]}.$$

\end{thm}
\begin{proof}
Let $\bar{U}$ denote the closure of $U$ in $X$, $\bar{j}: \bar{U} \hookrightarrow X$ the resulting closed immersion, and let $j': U \hookrightarrow \bar{U}$ denote the natural inclusion. Since $\bar{j}$ is a closed immersion, $\bar{j}_{!*} = \bar{j}_*$. Moreover, $j_{!*} = \bar{j}_{!*} \circ j'_{!*}.$ As a result, we may reduce to the case of an open immersion.\\

We may assume that $X$ (and $U$) is integral (i.e. connected and reduced). We may stratify $X$ by strata $S_i$ for $0 \leq i \leq d := \dim(X)$ such that:
\begin{enumerate}
    \item $\dim(S_i) = i$, each $S_i$ is smooth, and the closure $\bar{S_i} = \bigcup_{j \geq i} S_j$. 
    \item For each $-d \leq k \leq 0$, let $U_k := \bigcup_{i \leq k} S_{-i}$. We may find a stratification such that $U = U_r$ for some $r$. Note that $U_0 =X.$
\end{enumerate}
Let $j_{k-1}: U_{k-1} \hookrightarrow U_k$ denote the natural open immersions. Recall $U =U_r$ and $j: U \hookrightarrow X$ is the natural inclusion. With this notation, one has the following formula (see \cite{BBD}, 2.1.11):
$$j_{!*}(\s{G}) = \tau_{\leq -1}j_{-1*} \circ \tau_{\leq -2}j_{-2*} \circ \cdots \circ \tau_{\leq r}j_{r*}(\s{G}).$$
The result is now a direct consequence of Theorem \ref{thm:mainthm2}.

\end{proof}

We give an application of the previous result to the monodromy of intersection cohomology.
Let $X \rightarrow S$ be a proper morphism to a proper curve, $s_0 \in S$, and $j: U \hookrightarrow X$ be a smooth dense open subscheme. Let $\s{G} \in \s{P}(U)$, and consider $\s{H} := j_{!*}(\s{G}) \in \s{P}(X)$. Consider a loop $\gamma \in {\rm AL}(S)$ centered at $s_0$ and the corresponding map $h: \Delta \rightarrow S$. Up to shrinking the disk, we may assume that $Rf_*(\s{H})$ is locally constant when restricted to the punctured disk. For $t \in \Delta^*$, one has $R^if_*(\s{H})_t = \H^i(X_t, \s{H}_t)$, and the standard monodromy action of $\gamma$ on $\H^i(X_t, \s{H}_t)$. The previous theorem has the following corollary.

\begin{cor}\label{cor:intersectioncohom}
With notation as above, there is a finite set (denoted by $M$) of pairs $(\gamma', n_{\g'})$ with $\gamma' \in {\rm AL}(U)$, $n_{\gamma'}$ a positive integer and such that 
$${\rm Sp}(\gamma, \H^{i}(X_t,\s{H}_t)) \subset  \sum_{(\gamma',n_{\gamma'}) \in M} {\rm Sp}(\gamma',\s{G})^{[1/n_{\gamma'}]}.$$

\end{cor}

\section{Monodromy of Generalized Alexander Modules}\label{sec:alexmodules}
In this section, we explain how to deduce a local monodromy theorem in the setting of Alexander modules and discuss applications to computing monodromy of abelian covers. \\

\subsection{Monodromy of Alexander Modules}
Let $S$ be a smooth (connected) curve, and let $\pi: G \rightarrow S$ be a semi-abelian scheme. Consider a commutative diagram:
\begin{equation}\label{eqn:semiabeliansetting}
\xymatrix{
X \ar[dr]_{f} \ar[r]^{F}  & G \ar[d]^{\pi}\\
      & S
   }
\end{equation}
In this setting, one has the following data:
\begin{enumerate}
\item Let $e: S \rightarrow G$ denote the identity section. Consider the relative tangent bundle $\s{T}_{G/S}$, and vector bundle $e^*\s{T}_{G/S}$ on $S$. We have a commutative diagram:
$$
\xymatrix{
1 \ar[r] & \s{K} \ar[r] & e^*\s{T}_{G/S} \ar[rd] \ar[r]^{exp} & G \ar[r] \ar[d]& 1 \\
 &  &  & S  &  
}
$$
where the $exp$ is the exponential map, and $\s{K}$ is the kernel of the exponential map. By abuse of notation, we use the same notation $\s{K}$ to denote the sheaf of sections of $\s{K}$. This is a sheaf of abelian groups on $S$ with stalks $\s{K}_{s} = \pi_1(G_s,e(s))$ for a closed point $s \in S$. We set $\s{R}_S:= \Z[\s{K}]$. In particular, $\s{R}_{S,s} = \Z[\pi_1(G_s,e(s))]$.
\item Let $\s{L}_G:=(exp)_{!}(\Z)$. In the case where $S = \Sp(\C)$, $G$ is a semi-abelian variety, and $\s{L}_G$ is the the local system on $G$
whose stalk at $y \in G$ is given by the free abelian group on homotopy classes of paths from $e$ to $y$: 
$$(\s{L}_{G})_y = \Z[\pi_1(G;e,y)].$$
We view this as a (left) $R$-module, where $R = \Z[\pi_1(G,e)]$.
\item In general, for $s \in S$, we have $\s{L}|_{G_s} = \s{L}_{G_s}$. We also note that, by construction, $\s{L}_G$ is a sheaf of $\pi^{-1}(\s{R}_S)$-modules.

\item Below, we make the following additional hypothesis:
\begin{center}
(H) The semi-abelian scheme $G$ is an extension of an abelian scheme $A \rightarrow S$ by a torus $T \rightarrow S.$ We do not assume that $T$ is a split Torus.
\end{center}
It follows that $\pi$ is a fibre bundle. In particular, $\s{R}_S$ is locally constant. 
 
\item Let $s_0 \in S$, $\gamma \in {\rm AL}(S \setminus s_0)$ denote a (non-trivial) algebraic loop centered $s_0$, and $h: \Delta \rightarrow S$ denote corresponding map from the disk with center $h(0) 
 = s_0$. By abuse of notation, we use the same notation $h: \Delta^{*} \rightarrow S \setminus s_0$ to denote the restriction of $h$ to the corresponding punctured disk. By (4), the restriction of $G$ over the disk is a topological fibration. 

\item Consider now diagram \ref{eqn:semiabeliansetting} above but with everything restricted to $\Delta$. Then, under the hypothesis (H), $\s{R}_{\s{S}}$ can be (canonically) identified with the constant local system given by $R := \s{R}_{\s{S},s_0}$ (on the disk $\Delta$). Moreover, $\s{L}_{G}$ (restricted to $G_{\Delta}$) is a local system of $R$-modules. 
\item Let $\s{F}$ on $X$ be a constructible sheaf of $B$-modules where $B$ is a commutative noetherian ring of finite global dimension. We may consider 
the sheaf $\s{F}_R:= \s{F} \otimes_{\Z} F^{*}(\s{L}_{G})$. Under our hypothesis (H), and restricting to $\Delta^*$, this is a constructible sheaf of $B_R: = B \otimes_{\Z} R$-modules on $X$. 
\end{enumerate}

We wish to apply Theorem  \ref{thm:maintthm} and its variant Theorem \ref{thm:mainthmvariant}
to understand the monodromy action on 
$Rf_*\s{F}_R$. With $\g\in {\rm AL}(S)$ chosen above, we therefore consider $\g_X\in {\rm AL}(X)$ as in Theorem \ref{thm:mainthmvariant}.
Now, Definition 
\ref{defn:analyticloop} gives rise to the following three schemes:
\begin{enumerate}
 \item the closed subscheme ${\rm Sp}(\g_X,\s{F}) \subset \Sp(B[x])$,
 \item the closed subscheme $\Sp(\g_X,F^*\s{L}) \subset \Sp(R[x])$, and
 \item the closed subscheme $\Sp(\g_X,\s{F}_R) \subset \Sp(B_R[x]).$
\end{enumerate}

Recall, $B_R = B \otimes_{\Z} R$. The $B_R$-algebra homomorphism $B[x] \otimes_{\Z} R[x] \to B_R[x]$ given by $1  \otimes x\mapsto x$ and $x \otimes 1 \mapsto x$
induces $\t{diag}:\Sp(B_R[x]) \to \Sp(B[x]) \times \Sp(R[x])$. 

\begin{dfn}
Given closed subschemes $Z\subset\Sp\,B[x]$ and 
$W\subset \Sp\,R[x]$, 
we define $Z\overset{\bullet}{\times}W:=\t{diag}^{-1}(Z\times W)$. Furthermore, when $W=\Sp(R[x]/(x-M))$ for some $M\in R$, we will denote $Z\overset{\bullet}{\times}W$ by
$Z\overset{\bullet}{\times}M$. 
\end{dfn}

Note that $\s{F}_R=\s{F}\otimes_{\Z} F^*\s{L}$ implies that ${\rm Sp}(\g_X, \s{F}_R) ={\rm Sp}(\g_X,\s{F})\overset{\bullet}{\times}{\rm Sp}(\g_X,F^*\s{L}) $. The loop $\gamma_X$ maps to (via composition by $F$) a loop $\g_G\in {\rm AL}(G)$. It follows that ${\rm Sp}(\g_X,F^*\s{L})={\rm Sp}(\g_G,\s{L})$. The latter is determined by the homotopy class $[\g_G]\in\pi_1(\pi^{-1}\Delta)$.
Now, the loop $\gamma_G$ maps to a loop $\g_A\in {\rm AL}(A)$. The properness of $A\to S$
implies that $\g_A:\Delta^*\to A_{\Delta}$ extends to a map $\g_A:\Delta\to A_{\Delta}$, where $A_{\Delta}$ denotes the inverse image of $\Delta$ under $A\to S$.
It follows that the homotopy class $[\g_G]$
lies in the kernel of $\pi_1(\pi^{-1}\Delta)\to\pi_1(A_{\Delta})$ and the latter is clearly given by $\pi_1(T_{s_0})\hkr\pi_1(\pi^{-1}\Delta)$.
By the \emph{group of monomials} we mean the subgroup $\pi_1(T_{s_0})$ of the units of $R$. We now have $[\g_G]=M\in R\units$. In view of the fact that $\s{L}$ is a sheaf of free rank one $R$-modules, we see that
$\Sp(\g_G,\s{L})={\rm Sp}(R[x]/(x-M))$. By the above discussion, we have established that
$\Sp(\g_X,\s{F}_R)=\Sp(\g_X,\s{F})\overset{\bullet}{\times}M$.
We now apply Theorem \ref{thm:mainthmvariant} (1) to deduce:
\bt \label{thm:alexandermon}
With notation as above, given $\g\in AL(S)$ centered at $s_0\in S$, there are  
\begin{enumerate}
\item[(a)] loops $\g_i\in AL(X)$ and natural numbers $r_i$, for all $1\leq i\leq m$, and
\item[(b)] monomials $M_1,M_2,...,M_m\in \pi_1(T_{s_0})$
\end{enumerate}
such that the closed subscheme ${\rm Sp}(\g,R^qf_*\s{F}_R)$ of $\Sp(B_R[x])$
is contained in the sum of its closed subschemes $({\rm Sp}(\g_i,\s{F})\overset{\bullet}{\times}M_i)^{1/r_i}$ taken over $i=1,2,\ldots,m$:
$${\rm Sp}(\g,R^qf_*\s{F}_R) \subset \sum_{i=1}^m  ({\rm Sp}(\g_i,\s{F})\overset{\bullet}{\times}M_i)^{1/r_i}.$$
\et
\br
Suppose $B = K$, and $\s{F}$ is quasi-unipotent. If $G =A$, then the group of monomials is trivial, and it follows from the previous corollary that the monodromy action on $R^qf_*\s{F}_R$ is quasi-unipotent i.e. the eigenvalues of the monodromy action are roots of unity. 
\er

\bex\label{example:alexmod}
We may apply the previous theorem to the following geometric setting. Let $Y \subset X$ be a closed subvariety and consider $j: X \setminus Y \hookrightarrow X$. We set $\s{F} := j_!j^*\Z$, and let $\s{F}_R$ be as in the Theorem above. With notation as above, we have a local system $R^if_*(\s{F}_R)$ of $R$-modules (after restriction to a sufficiently small disk). If $X \rightarrow S$ is proper, or if over disk we have $X$ is a topological fibration, then
$R^if_*(\s{F}_R)_t = \H^i(X_t,Y_t; F^*_t\Z[\pi_1(G_t,e(t)])$ for a general $t \in \Delta^*$ and the corresponding monodromy representation. The above theorem reduces us to computing the monodromy of the corresponding universal local system $F^{-1}(\s{L}_{G})$, and therefore of $\s{L}_{G}$. In particular, we find that the eigenvalues of monodromy are given by $M \in R$ whose $r$-th power is a monomial (from Torus). In particular, we obtain results for the monodromy action on `generalized Alexander modules'. If $G = A$, then by the previous remark we get roots of unity, i.e. it is quasi-unipotent.  
\eex

\subsection{Abelian Coverings}
We continue with the notation and hypotheses of the previous section. In particular, $X, S, G,$ and $\s{F}$ etc. are as in the previous section. We fix a loop $\gamma \in {\rm AL}(S)$ centered at $s_0$, and work over a disk $\Delta$ as before.\\

Given a finite etale morphism $\phi:H\to G$, its base change $F^*\phi:X\times_GH\to X$ is also a finite etale morphism. Let $n_G:G\to G$ denote multiplication by a natural number $n$; we denote the base change $F^*n_G$ by $n_X:X_n\to X$ and let $f_n:X_n\to S$ denote the resulting composition given by $f_n=f\circ n_X$. Finally, we consider the sheaf 
$\s{F}_n=n_X^*\s{F}$ on $X_n$. We have the resulting commutative diagram:

$$
\xymatrix{
X_n \ar[d]_{n_X}  \ar[r]^{F_n} &  G \ar[d]^{n_G} \\
X   \ar[r]^{F} \ar[dr]_{f} &   G \ar[d]^{\pi}   \\
    &   S.}
$$

\begin{question}
 What is the {\it local monodromy} of $R^q(f_n)_*\s{F}_n$ at the loop $\gamma$? If $\s{F}$ is the constant local system $\Z$, then these are roots of unity. What roots of unity appear?
\end{question}

This is essentially the question \ref{question:questionintro} stated in the introduction. We begin by making the question above more precise.
First, note that $R^q(f_n)_*\s{F}_n=R^qf_*({n_X}_*\s{F}_n)$ since $n_X$ is a finite morphism. By the projection formula
$${n_X}_*\s{F}_n=\s{F}\otimes {n_X}_*\Z_{X_n}.$$
On the other hand (working over $\Delta$),
$$\s{F}\otimes {n_X}_*\Z_{X_n} = \s{F}\otimes F^*R_n $$ where $R_n$ is the group-ring 
of $V/V^n$ where $V$ is the fundamental group of $\pi^{-1}s_0$. Note that $R_n$ is viewed as a local system on $G$. With this notation, Theorem \ref{thm:alexandermon} now has the following corollary:

\begin{cor}\label{cor:abcovers}
 With notation as above, $M:=R^q(f_n)_*\s{F}_n$ has the natural structure of a $B_n:=B\otimes R_n$-module, and the $\g$-action is an automorphism of this module. Let ${\rm Sp}(\gamma, \s{F}_n)$ denote the corresponding closed subscheme of $\Sp(B_n[x])$. Then there are 
 \begin{enumerate}
\item[(a)] loops $\g_i\in {\rm AL}(X)$, and natural numbers $r_i$, for all $1\leq i\leq m$, and
\item[(b)] monomials $M_1,M_2,...,M_m\in \pi_1(T_{s_0})$
\end{enumerate}
such that 

$${\rm Sp}(\g,R^qf_*\s{F}_R) \subset \sum_{i=1}^m  ({\rm Sp}(\g_i,\s{F})\overset{\bullet}{\times}M_i)^{1/r_i}.$$
and 
$${\rm Sp}(\g,R^q(f_n)_*\s{F}_n) \subset \sum_{i=1}^m  ({\rm Sp}(\g_i,\s{F})\overset{\bullet}{\times}M_{i,n})^{1/r_i},$$
$M_{i,n}$ is the image of $M_i$  in $R_n\units$. 
\end{cor}
\begin{proof}
This follows immediately from Theorem \ref{thm:mainthmvariant}.
\end{proof}

We now explain how to use the corollary above in order to solve the question of the introduction. In that case, we take $\s{F} = \Z$ as the constant local system. By the local monodromy theorem, we know that the eigenvalues of the local monodromy of $R{f_n}_*\Z_{X_n}$ are roots of unity. The above theorem helps to answer which roots of unity appear. More precisely, we look at the sheaf $\s{F}_R$, compute the corresponding monomials $M$ and consider their images $M_n$. In particular, this gives a {\it uniform in $n$} computation of the eigenvalues of local monodromy.

\bibliographystyle{plain}
\bibliography{nonabelian}

\end{document}